\documentclass[11pt]{article}

\usepackage[T1]{fontenc}
\usepackage[utf8]{inputenc}
\usepackage{amsmath,amssymb,epsfig,bbm}
\usepackage{stmaryrd}
\usepackage[colorlinks]{hyperref}
\usepackage{comment}
\usepackage{color}
\usepackage{textcomp}

\usepackage[acronym,toc]{glossaries}
\usepackage{glossary-longragged}
\usepackage{glossary-superragged}
\newglossary[slg]{symbols}{sym}{sbl}{List of Symbols}

\usepackage[textsize=small]{todonotes}
\usepackage{varwidth}

\usepackage{url}

\usepackage{enumitem}
\setlist[itemize]{labelindent=*, leftmargin=.5 truecm,nosep}

\usepackage[nottoc]{tocbibind}

\usepackage{amsthm}

\theoremstyle{plain}
\newtheorem{defi}{Definition}[section]
\newtheorem{prop}[defi]{Proposition}
\newtheorem{theo}[defi]{Theorem}
\newtheorem{conj}[defi]{Conjecture}
\newtheorem{lemm}[defi]{Lemma}
\newtheorem{coro}[defi]{Corollary}

\theoremstyle{definition}
\newtheorem{rema}[defi]{Remark}
\newtheorem{exem}[defi]{Example}
\newtheorem{exems}[defi]{Examples}

\newcommand{\bdefi}{\begin{defi}}
\newcommand{\edefi}{\end{defi}}
\newcommand{\bprop}{\begin{prop}}
\newcommand{\eprop}{\end{prop}}
\newcommand{\btheo}{\begin{theo}}
\newcommand{\etheo}{\end{theo}}
\newcommand{\blemm}{\begin{lemm}}
\newcommand{\brema}{\begin{rema}}
\newcommand{\erema}{\end{rema}}
\newcommand{\bexer}{\begin{exem}}
\newcommand{\eexer}{\end{exem}}
\newcommand{\bexem}{\begin{exem}}
\newcommand{\eexem}{\end{exem}}
\newcommand{\bexems}{\begin{exems}}
\newcommand{\eexems}{\end{exems}}
\newcommand{\bconj}{\begin{conj}}
\newcommand{\econj}{\end{conj}}
\newcommand{\elemm}{\end{lemm}}
\newcommand{\bcoro}{\begin{coro}}
\newcommand{\ecoro}{\end{coro}}
\newcommand{\dem}{\noindent{\bf Proof. }}
\newcommand{\rem}{\noindent{\bf Remark. }}

\usepackage{mathrsfs}
\renewcommand\mathcal{\mathscr}

\newcommand{\C}{{\cal C}}
\newcommand{\D}{{\cal D}}

\newcommand{\F}{{\cal F}}
\newcommand{\G}{{\cal G}}
\renewcommand{\H}{{\cal H}}
\newcommand{\I}{{\cal I}}

\newcommand{\M}{{\cal M}}
\newcommand{\N}{{\cal N}}
\newcommand{\OOO}{{\cal O}}

\newcommand{\maths}[1]{{\mathbb #1}}  

\newcommand{\CC}{\maths{C}}
\newcommand{\DD}{\maths{D}}
\newcommand{\FF}{\maths{F}}
\newcommand{\HH}{\maths{H}}

\newcommand{\NN}{\maths{N}}

\newcommand{\PP}{\maths{P}}
\newcommand{\QQ}{\maths{Q}}
\newcommand{\RR}{\maths{R}}
\newcommand{\SSS}{\maths{S}}

\newcommand{\UU}{\maths{U}}
\newcommand{\XX}{\maths{X}}

\newcommand{\ZZ}{\maths{Z}}

\renewcommand{\ggg}{{\mathfrak g}}

\newcommand{\mmm}{{\mathfrak m}}

\newcommand{\bs}{\backslash}
\newcommand{\ga}{\gamma}
\newcommand{\Ga}{\Gamma}
\newcommand{\ov}[1]{{\overline{#1}}} 
\newcommand{\ra}{\rightarrow}
\newcommand{\wt}[1]{{\widetilde{#1}}}
\newcommand{\wh}[1]{{\widehat{#1}}}

\newcommand{\Aut}{\operatorname{Aut}}

\newcommand{\bigO}{\operatorname{O}}
\newcommand{\card}{{\operatorname{Card}}}

\newcommand{\cqfd}{\hfill$\Box$}

\newcommand{\Dirac}{{\Delta}}

\newcommand{\gengeod}
{\operatorname{\widecheck{\G\,}\!\!}}

\newcommand{\id}{\operatorname{id}}
\renewcommand{\Im}{{\operatorname{Im}}}
\newcommand{\Isom}{\operatorname{Isom}}

\renewcommand{\log}{\operatorname{ln}}

\newcommand{\n}{\operatorname{\tt n}}

\renewcommand{\Re}{{\operatorname{Re}}}

\newcommand{\semidirect}{{\mathbb n}}

\newcommand{\tr}{\operatorname{\tt tr}}

\newcommand{\Vol}{\operatorname{Vol}}
\newcommand{\vol}{\operatorname{vol}}

\newcommand{\hdr}{{\HH}^2_\RR}

\newcommand{\htr}{{\HH}^3_\RR}

\newcommand{\hnr}{{\HH}^n_\RR}

\newcommand{\hdc}{{\HH}^2_\CC}

\newcommand{\hdh}{{\bf H}^2_\HH}

\newcommand{\GL}{\operatorname{GL}}
\newcommand{\PGL}{\operatorname{PGL}}
\newcommand{\SL}{\operatorname{SL}}
\newcommand{\PSL}{\operatorname{PSL}}

\newcommand{\SU}{\operatorname{SU}}
\newcommand{\PU}{\operatorname{PU}}
\newcommand{\PSU}{\operatorname{PSU}}

\newcommand{\SLC}{\operatorname{SL}_{2}(\CC)}
\newcommand{\PSLC}{\operatorname{PSL}_{2}(\CC)}

\newcommand{\PSLR}{\operatorname{PSL}_{2}(\RR)}

\newcommand{\PSLZ}{\operatorname{PSL}_{2}(\ZZ)}

\newcommand{\flow}[1]{{{\tt g}^{#1}}}  
\newcommand{\wss}{W^+} 
\newcommand{\wsu}{W^-} 
\newcommand{\ws}{W^{0+}} 
\newcommand{\wu}{W^{0-}} 
\newcommand{\muss}[1]{{\mu_{\wss(#1)}}}
\newcommand{\mus}[1]{{\mu_{\ws(#1)}}}

\newcommand\normalout{\partial^1_{+}}
\newcommand\normalin{\partial^1_{-}}
\newcommand\normalpm{\partial^1_{\pm}}
\newcommand\normalmp{\partial^1_{\mp}}

\usepackage{interval}

\usepackage{mathtools}
\makeatletter
\DeclareRobustCommand\widecheck[1]{{\mathpalette\@widecheck{#1}}}
\def\@widecheck#1#2{%
    \setbox\z@\hbox{\m@th$#1#2$}%
    \setbox\tw@\hbox{\m@th$#1%
       \widehat{%
          \vrule\@width\z@\@height\ht\z@
          \vrule\@height\z@\@width\wd\z@}$}%
    \dp\tw@-\ht\z@
    \@tempdima\ht\z@ \advance\@tempdima2\ht\tw@ \divide\@tempdima\thr@@
    \setbox\tw@\hbox{%
       \raise\@tempdima\hbox{\scalebox{1}[-1]{\lower\@tempdima\box
\tw@}}}%
    {\ooalign{\box\tw@ \cr \box\z@}}}
\makeatother

\newcounter{fig}

\title{Joint partial equidistribution of Farey rays\\
  in negatively curved manifolds and trees}
\author{Jouni Parkkonen \and Fr\'ed\'eric Paulin} 
\date{\today}

\begin{document}
\bibliographystyle{../alphanum}
\maketitle
\begin{abstract}
  We prove a joint partial equidistribution result for common
  perpendiculars with given density on equidistributing equidistant
  hypersurfaces, towards a measure supported on truncated stable
  leaves. We recover a result of Marklof on the joint partial
  equidistribution of Farey fractions at a given density, and give
  several analogous arithmetic applications, including in Bruhat-Tits trees.
\end{abstract}

\section{Introduction}
\label{sec:intro}
In this paper, we study geometric equidistribution results on
negatively curved manifolds with applications to arithmetic problems.
Let $M$ be a complete connected Riemannian manifold with pinched
negative sectional curvature at most $-1$. Let $m_{\rm BM}$ be its
Bowen-Margulis measure, which, when finite and renormalized and when
the sectional curvature has bounded derivative, is the probability
measure of maximal entropy for the geodesic flow on $T^1M$.  For
example when $M$ has finite volume, it is well-known that the
conditional measure of $m_{\rm BM}$ on the image $\flow{t} W$ of a
closed strong unstable leaf $W$ by the geodesic flow $\flow{t}$ at
time $t$ equidistributes towards $m_{\rm BM}$ as $t\ra +\infty$.  See,
for instance, the works of Dani, Eskin-McMullen
\cite[Thm.~7.1]{EskMcMul93}, Margulis, Kleinbock-Margulis
\cite[Prop.~2.2.1]{KleMar96}, Ratner, Sarnak \cite[Thm.~1]{Sarnak81},
as well as \cite[Thm.~1]{ParPau14ETDS} and
\cite[Thm.~10.2]{BroParPau19} for generalisations.  Given a natural
family $(\F_t)_{t\in\RR}$ of finite subsets $\F_t$ of points on
$\flow{t} W$, a natural question is to study the limiting distribution
properties of $\F_t$ as $t\ra+\infty$. If $\F_t$ is denser and denser
in $\flow{t} W$, it is expected that $\F_t$ will also equidistribute
to $m_{\rm BM}$. If $\F_t$ is too sparse in $\flow{t} W$, it is
expected for the limiting distribution to be purely punctual. A
threshold seems to occur when $\F_t$ has a constant density in
$\flow{t} W$, possibly yielding equidistribution of partial nature.

In this paper, we take $\F_t$ to be the image by $\flow{t}$ of the
subset of $W$ of initial tangent vectors of the common perpendiculars
to another cusp neighbourhood, having a length bound chosen in order
to have a constant density at each time $t$. We prove that $\F_t$ then
equidistributes towards the conditional measure of $m_{\rm BM}$ on a
truncated weak stable leaf. We, for instance, recover the case $n=2$
of a theorem by Marklof \cite[Thm.~6]{Marklof10Inv}, as well as
\cite[Thm.~6.1]{Lutsko22}. We actually prove a joint partial
equidistribution result, for more general families, give a
version of our results for tree quotients, and give several arithmetic
applications.

More precisely, let $\wt M$ be a complete simply connected Riemannian
manifold with pinched negative sectional curvature at most $-1$, and
let $\Ga$ be a nonelementary discrete subgroup of $\Isom(X)$, with
critical exponent $\delta_\Ga$ (see for instance \cite{BriHae99}). Let
$D$ be a nonempty proper closed convex subsets of $\wt M$ and $H$ an
horoball of $\wt M$ such that the families $\D^-=(\ga D)_{\ga\in\Ga}$
and $\D^+=(\ga H)_{\ga\in\Ga}$ are locally finite in $\wt M$.

Let us introduce the measures that come into play in this paper,
refering to Section \ref{sec:background} and \cite{BroParPau19} for
more explanations. We denote by $\|\mu\|$ the total mass of a measure
$\mu$.

Let $(\mu_{x})_{x\in \wt M}$ be a Patterson density for $\Ga$ and let
$m_{\rm BM}$ be the associated Bowen-Margulis measure on $\Ga\bs
T^1\wt M$.  When $\wt M$ is a symmetric space and $\Ga$ has finite
covolume, then $\mu_{x}$ is (a scalar multiple of) the unique
probability measure on $\partial_\infty \wt M$ invariant under the
stabiliser of $x$ in the isometry group of $\wt M\,$, and $m_{\rm BM}$
is the Liouville measure, which is finite and mixing. Let $W$ be the
strong stable leaf in $T^1\wt M$ whose image in $\wt M$ is $\partial
H$, and let $\mu^{0+}_{\D^+,t_0}$ be the conditional measure of
$m_{\rm BM}$ on the truncated weak stable leaf $\Ga \bigcup_{s\geq
  t_0} \flow{s}W$.  The measure $\mu^{0+}_{\D^+,t_0}$ is finite and
nonzero for instance when $H$ is centered at a bounded parabolic fixed
point of $\Ga$. Let $\sigma^+_{\D^-}$ be the outer skinning measure of
$\D^-$, see for instance \cite{ParPau14ETDS}, as well as
\cite{OhSha12,OhSha13} when $\wt M$ is geometrically finite with
constant curvature, and when $D$ is a ball, horoball or complete
totally geodesic submanifold.  When $D$ is a horoball,
$\sigma^+_{\D^-}$ is the conditional measure of $m_{\rm BM}$ on the
strong unstable leaf in $T^1M$ having a lift to $T^1\wt M$ whose image
in $\wt M$ is $\partial D$.

For every $\ga\in\Ga$ such that $d(D,\ga H)>0$, let $v_\ga\in T^1\wt
M$ be the normal outgoing vector of $D$ pointing towards the point at
infinity of $\ga H$.
    
\btheo\label{theo:mainintro} Assume that $m_{\rm BM}$ is finite and
mixing for the geodesic flow on $\Ga\bs T^1\wt M$, and that
$\sigma^{+}_{\D^-}$ and $\mu^{0+}_{\D^+, t_0}$ are finite and nonzero.
Then for every $t_0\in\RR$, for the weak-star convergence of measures
on $(\Ga\bs T^1\wt M)^2$, we have
$$
\lim_{t\ra+\infty} \;\|m_{\rm BM}\|\;e^{-\delta_\Ga\, t}
\sum_{\substack{\ga\in\Ga_{D}\!\backslash \Ga/\Ga_{H}\\ 
0<d(D,\ga H)\leq t-t_0}} \; 
\Delta_{\Ga v_\ga}\otimes\Delta_{\flow t\Ga v_\ga}
\;=\sigma^+_{\D^-}\otimes\mu^{0+}_{\D^+,t_0}\;.
$$
\etheo

See Theorem \ref{theo:quotient} for a more general version (with
improvement to narrow convergence under additional hypotheses), as
well as a version for quotients of trees by discrete groups of
automorphisms. See Section \ref{sec:shrinking} for a proof, after some
preliminary work in Section \ref{sec:background}, in particular on the
truncated weak stable leaves and their measures. The proof starts by
using the joint equidistribution result of common perpendiculars from
\cite{ParPau17ETDS}, but the statement of Theorem 1.1 is only
seemingly similar to Eq.~(12) in loc.~cit.~and new ideas and techniques 
are required. One of these ideas is a new subdivision scheme along the
geodesic flow that allows a good control of the exponential growth. One
of the techniques is an important regularity study of the splitting of
the weak stable leaves and of the dynamics on the unstable
horospheres.
    
As a consequence of our main result Theorem \ref{theo:mainintro}, we
recover the case $n=2$ of a theorem by Marklof
\cite[Thm.~6]{Marklof10Inv} on the joint partial equidistribution of
Farey points chosen with constant average density on an
equidistributing horocycle on the modular curve $\PSLZ\bs\hdr$, see
Corollary \ref{coro:marklof}.  In the present case (contrarily to
other distribution results in number theory), the restriction to a
fixed denominator of the Farey fractions in \cite{Marklof10Inv} is
only marginally stronger, by the growth properties of the
horospheres. The relationship between Farey fractions and hyperbolic
geometry (and in particular with the divergent geodesics) is not new,
probably going back to Ford. See for instance the works of
Athreya-Cheung \cite{AthChe14}, Sarnak, Series, Sullivan, and the
references of \cite{HerPau02b,ParPauCIRM}.  We also recover
\cite[Thm.~6.1]{Lutsko22}, originally proved for hyperbolic
surfaces. 

We give in Section \ref{sec:applicationexamples} several
generalisations of Marklof's result, as for instance the complex
version below. See Corollary \ref{coro:complexFarey} for a more
general statement, and Subsections \ref{subsec:heisenbergFarey} and
\ref{subsec:heisenbergFareyquat} for distribution results of Farey
points with constant average density on closed horospheres in complex
and quaternionic hyperbolic orbifolds.  It might be that it is
possible to obtain these applications using purely homogeneous
dynamics techniques, along the lines of the cross-section method of
Marklof \cite{Marklof10Inv} and Athreya-Cheung \cite{AthChe14}. But no
such results appear in the literature yet. We believe that covering
all our examples might require a lot of work, even starting from the
complex case with a large class number of the imaginary quadratic
field, as the cross-sections, as well as other fundamental domain
issues, are highly more complicated for general arithmetic lattices in
rank one real Lie groups than for $\SL_2(\ZZ)$. Furthermore, the case
of groups over local fields with positive characteristic is likely to
require major innovations by homogeneous dynamics methods.

Let $K$ be an imaginary quadratic number field, with ring of integers
$\OOO_K$ and discriminant different from $-4$ and $-3$ in order to
simplify the statement in this introduction.  Let $G=\PSLC$, let $\Ga$
be the Bianchi group $\PSL(\OOO_K)$, let
\[
H=\big\{n_-(r)=\begin{bmatrix} 1 & r\\0&1\end{bmatrix}:r\in\CC\big\}\quad 
{\rm and}
\quad 
\forall\;t\in \RR,\;\Phi^t=
\begin{bmatrix} \;e^{-t/2} & 0\\0&e^{t/2}\end{bmatrix}\,.
\]
Let $M=\Big\{\begin{bmatrix} \;e^{-i\,\theta/2} &
0\\0&e^{i\,\theta/2}\end{bmatrix}:\theta\in\RR\Big\}$.  We endow the
compact abelian groups $\CC/\OOO_K$ and $(H\cap\Ga)\bs H$ with their
probability Haar measures $dx$ and $d\mu_{(H\cap\Ga)\bs H}$.  For every
$t\in\RR$, we consider the set $\F_{t}$ of {\em complex Farey
  fractions of height at most} $e^{t/2}$, defined by
\[
\F_{t}=\Big\{\frac pq \mod \OOO_K:\;p,q\in\OOO_K,
\quad p\OOO_K+q\OOO_K=\OOO_K,\quad 0<|q|\leq e^{t/2}\Big\}\,.
\]

\bcoro\label{coro:complexFareyintro}
Let $f:(\CC/\OOO_K) \times(\Ga\bs G/M)\to\RR$ be a bounded
continuous function. Then for every $t_0\in\RR$, we have
\begin{align*}
\lim_{t\to+\infty}\;\;&\frac 1{\card\,\F_{t-t_0}}
\sum_{r\in\F_{t-t_0}}f(r,\Ga n_-(r)\Phi^tM)\\&=2\,e^{2\,t_0}
\int_{s=t_0}^{+\infty}\int_{y\,\in\,(H\cap \Ga)\bs H} \int_{x\,\in\,\CC/\OOO_K} 
f(x,\;^ty^{-1}\Phi^{s}M) \,dx\,d\mu_{(H\cap \Ga)\bs H}(y)\,e^{-2\,s}\,ds\,.
\end{align*}
\ecoro

We now give a joint partial equidistribution result of arithmetic
points with given density on an expanding horosphere in an arithmetic
quotient of a nonarchimedean simple Lie group (see Corollary
\ref{coro:nonarchFarey} for a more general version). Let $R=\FF_q[Y]$
be the ring of polynomials over a finite field $\FF_q$ with one
indeterminate $Y$, and let $\wh K=\FF_q((Y^{-1}))$ be the valued field
of formal Laurent series in $Y^{-1}$ over $\FF_q$ with
$|Y^{-1}|=1$. Let $G=\PGL_2(\wh K)$, let $\Ga=\PGL(R)$, let
\[
H=\Big\{n_-(r)=\begin{bmatrix} 1 & r\\0&1\end{bmatrix}:r\in\wh K\Big\}\quad 
{\rm and}
\quad 
\forall\;n\in \ZZ,\;\Phi^n=
\begin{bmatrix} \;1 & 0\\0&Y^{n}\end{bmatrix}\,.
\]
Let $\Ga_H=N_G(H)\cap \Ga$ and $M=\Big\{\begin{bmatrix} \;1 &
0\\0&u\end{bmatrix}:u\in\wh K,|u|=1\Big\}$.  We endow $\Ga_H\bs H$
with the induced measure $d\mu_{\Ga_H\bs H}$ of a Haar measure of $H$,
normalised to be a probability measure.  For every $n\in\ZZ$, we
consider the set $\F_{n}$ of {\em nonarchimedean Farey fractions of
  height at most} $q^n$, defined by
\[
\F_{n}=\Ga_H\bs \Big\{n_-\big(\frac PQ\big):\;P,Q\in R,
\quad PR+QR=R,\quad 0\leq\deg Q\leq n\Big\}\,.
\]

\bcoro\label{coro:nonarchFareyintro}
Let $f:(\Ga_H\bs H) \times(\Ga\bs G/M)\to\RR$ be a bounded
continuous function. Then for every $n_0\in\ZZ$, we have
\begin{align*}
\lim_{n\to+\infty}\;\;&\frac 1{\card\,\F_{n-n_0}}
\sum_{r\in\F_{n-n_0}}f(r,\Ga\, r\,\Phi^{2n}M)\\
&=(q^2-1)\,q^{2n_0-2}
\sum_{m=n_0}^{+\infty}\int_{x,y\,\in\,\Ga_H\bs H} 
f(x,\Ga\;^ty^{-1}\Phi^{2m}M)
\,d\mu_{\Ga_H\bs H}(x)\,d\mu_{\Ga_H\bs H}(y)\,q^{-2m}\,.
\end{align*}
\ecoro

\medskip
\noindent{\small {\it Acknowledgements: } The authors thank for its
  support the French-Finnish CNRS IEA PaCap.}

\section{Background and definitions}
\label{sec:background}

Let $X$ be either a complete simply connected Riemannian manifold with
pinched negative sectional curvature at most $-1$ or a proper
geodesically complete $\RR$-tree.  Let $\Ga$ be a nonelementary
discrete group of isometries of $X$. We refer to \cite[Chap.~2 and
  3]{BroParPau19}, with potential $0$ throughout this paper, for
background informations on the data $(X,\Ga)$, and in particular for
the definition of the boundary at infinity $\partial_\infty X$ of $X$,
and of the critical exponent $\delta_{\Ga}>0$ of $\Ga$ in its Section
3.3.

We refer to \cite[\S2.2]{BroParPau19} for the following definitions.
We denote by $\gengeod X$ the Bartels-Lück space of generalised
geodesics in $X$ (that is, of continuous maps $\RR\ra X$ that are
isometric on a closed interval of $\RR$ with nonempty interior and
locally constant outside it), endowed with the distance $d$ defined by
\begin{equation}\label{eq:defdBartelLuck}
\forall\; \ell,\ell'\in \gengeod X,\;\;\;
d(\ell,\ell')=\int_{-\infty}^{+\infty}d(\ell(t),\ell'(t))\;e^{-2|t|}\;dt\;.
\end{equation}
It contains the closed subspace $\G X$ of (true) geodesic lines and
the closed subspaces $\G_{\pm,0} X$ of (positive/negative) geodesic
rays, that is, of generalised geodesics that are isometric on exactly
$\pm[0, +\infty[$ (that we identify with their restriction to
$\pm[0, +\infty[$).   We denote by $\ell\mapsto\ell_\pm$ the two
endpoint maps from $\gengeod X$ to $X\cup\partial_\infty X$. Let
$(\flow t)_{t\in\RR}$ be the (continuous-time) geodesic flow on
$\gengeod X$, which preserves $\G X$. Let $\G_\pm X=\G X\cup
\bigcup_{t\in\RR} \flow t \G_{\pm,0}$ be the closed subspace of
generalised geodesics that are isometric at least on an interval
$\pm[a,+\infty[$ for some $a\in\RR$, so that $\G_-X\cap \G_+X=\G X$.
The Bartels-Lück space is important in order to allow the positive
geodesic rays pushed by the geodesic flow at large positive times
to converge to geodesic lines.

We denote by $\wt m_{\rm BM}$ the Bowen-Margulis measure of $\Ga$ on
$\G X$ and by $m_{\rm BM}$ the Bowen-Margulis measure on $\Ga\bs\G X$
associated with any choice of Patterson-Sullivan density
$(\mu_x)_{x\in X}$, see for instance \cite{Roblin03} or
\cite[\S4.2]{BroParPau19} with potential $0$.

Given a proper closed convex subset $D$ of $X$, we refer to
\cite[\S2.4]{BroParPau19} for the definition of their inner/outer
normal bundles $\normalpm D$, which are contained in $\G_{\pm,0}
X$. We refer to \cite[\S7.1]{BroParPau19} again with potential $0$
(see also \cite{ParPau14ETDS} in the manifold case) for the definition
of the outer/inner skinning measures $\wt \sigma^\pm_D$ on $\normalpm
D$. Given a measurable map $f$, we denote by $f_*$ the pushforward map of
measures. Recall that, for every $\ga\in\Ga$, we have
\begin{equation}\label{eq:equivskinn}
 \ga_*(\wt \sigma^\pm_D)=\wt \sigma^\pm_{\ga D}\,.
\end{equation}

Given $w$ in $\G_+ X$ or $\G_-X$ respectively, we refer to
\cite[\S2.3]{BroParPau19} for the definitions of its strong stable
leaf $\wss(w)$ or strong unstable leaf $\wsu(w)$, of its (weak) stable
leaf $\ws (w)$ or (weak) unstable leaf $\wu (w)$, and of its stable
horoball $H\!B_+(w)$ or unstable horoball $H\!B_-(w)$. The {\em
  antipodal (or time reversal) map} $\iota:\gengeod X\ra\gengeod X$
defined by $\ell\mapsto \{t\mapsto \ell(-t)\}$ is an involution
satisfying $\iota(\G_+ X)=\G_- X$ and
\[
\forall\;w\in\G_+X, \quad \iota\,\wss(w)=\wsu(\iota\, w)\;.
\]

Let $w\in\G_+ X$. We refer to \cite[\S2.4]{BroParPau19} for the
definition of the canonical homeomorphism $N^+_w: \wss(w)\to\normalin
H\!B_+(w)$ that associates to a geodesic line $\ell\in\wss(w)$ the
unique (negative) geodesic ray $\rho \in \normalin H\!B_+(w)$ such
that $\ell_-=\rho_-$. We also denote by abuse $N^+_w(\ell)=\ell_{\mid
  \,]-\infty,0]}$. The homeomorphism $N^+_w$ relates the inner
skinning measure $\wt \sigma^-_{H\!B_{+}(w)}$ of $H\!B_+(w)$ to the
conditional $\muss{w}$ on the strong stable leaf $\wss(w)$ of $w$ of
the Bowen-Margulis measure $\wt m_{\rm BM}$ as follows: for $\ell\in
\wss(w)$, we have
\begin{equation}\label{eq:strongstabmeas}
d\muss{w}(\ell)=d\,((N^+_w)^{-1})_*\wt \sigma^-_{H\!B_{+}(w)}(\ell)=
d\,\wt \sigma^-_{H\!B_{+}(w)}(\ell_{\mid \,]-\infty,0]})\,.
\end{equation}

Recall that we have a homeomorphism 
$$h_w:\wss(w)\times\RR\ra \ws(w)
,\hspace{3mm}
(\ell,s)\mapsto \flow s\ell\,.
$$
For every isometry $\ga$ of
$X$, for all $t,s\in\RR$ and $\ell\in\wss(w)$, we have
$$
\ga h_w(\ell,s)=h_{\ga w}(\ga \ell,s)\;\;\;{\rm and}\;\;\;
  \flow t\circ h_w(\ell,s)=h_{\flow t w}(\flow t \ell,s)\;.
$$
  
The homeomorphism $h_w$ writes the conditional measure $\mus{w}$ on
the stable leaf $\ws(w)$ of $w$ of the Bowen-Margulis measure $\wt
m_{\rm BM}$ as a twisted product measure of the measure $\muss{w}$ on
$\wss(w)$ and the Lebesgue measure on $\RR$, see
\cite[Eq.~(7.12)]{BroParPau19} with potential $0$: for all $s\in\RR$
and $\ell\in \wss(w)$, we have
\begin{equation}\label{eq:stabmeas}
d\mus{w}(\flow s\ell)=e^{-\delta_{\Ga} s}d\muss{w}(\ell)\,ds\;.
\end{equation}
Note that for every $\ga\in\Ga$, we have
\begin{equation}\label{eq:equivweakstabmeas}
\ga_*\mus{w}=\mus{\ga w}\,.
\end{equation}
Since the Lebesgue measure is atomless, for every Borel subset
$\Omega^+$ of $\wss(w)$, the boundary of $h_w(\Omega^+\times [a,b])$
has measure $0$ for $\mus{w}$ if and only if the boundary of
$\Omega^+$ has measure $0$ for $\muss{w}$.

\medskip
For all $w\in\G_+ X$ and  $s\in\RR$, let
$$
\flow s_\mid :\normalin H\!B_+(w)\ra
\normalin H\!B_+(\flow s w)
$$
be the homeomorphism that associates to $\rho\in\normalin H\!B_+(w)$
the unique $\rho' \in \normalin H\!B_+(\flow s w)$ such that
$\rho_-=\rho'_-$, or equivalently such that we have $\rho(t)=
\rho'(t-s)$ for every $t\in\RR$ such that $t\leq 0$ and $t-s\leq 0$.
Note that $\flow s \wss(w)=\wss(\flow s w)$ and that the following
diagram is commutative
\begin{equation}\label{eq:diagcomNplus}
\begin{array}{ccc}
  \wss(w)&\stackrel{N^+_w}{\longrightarrow} &\normalin H\!B_+(w)\smallskip\\
  ^{\flow s} \downarrow & & \downarrow\;^{\flow s_\mid}\\
  \wss(\flow sw)&\stackrel{N^+_{\flow s w}}{\longrightarrow} &
  \normalin H\!B_+(\flow s w)\;.
\end{array}
\end{equation}

Let us now introduce the truncated (weak) stable leaves in $T^1\wt M$.
The projections on the second factor of the limiting measures of our
empirical joint distributions will have as support the union of a
locally finite family of truncated stable leaves. For every $\sigma\in
\RR\cup\{-\infty\}$, the {\it $\sigma$-stable leaf} of $w\in\G_+ X$ is
$$
\ws_\sigma(w)=\bigcup_{t\ge\sigma}\flow t \wss(w)\,,
$$
so that $\ws_{-\infty}(w)$ equals $\ws(w)$.

\blemm\label{lem:unifcontwss} Let $w\in \G_+ X$ and $s\in\RR$.
\begin{enumerate}
\item[(1)] The homeomorphism $N^+_{\flow{s}w}: \wss(\flow{s}w)\to
  \normalin H\!B_+(\flow{s}w)$ is uniformly bicontinuous, uniformly in
  $s$.
\item[(2)] The homeomorphism $h_w:\wss(w)\times\RR\ra \ws(w)$ is
  uniformly bicontinuous.
\item[(3)] The homeomorphism $\flow s_\mid :\normalin H\!B_+(w)\ra
  \normalin H\!B_+(\flow s w)$ is uniformly bicontinuous, uniformly on
  $s$ varying in a compact subset of $\RR$.
\end{enumerate}
\elemm

\dem (1) For all $\ell,\ell'\in \wss(\flow{s}w)$, by Equation
\eqref{eq:defdBartelLuck}, we have
\begin{align*}
  d(N^+_{\flow{s}w}(\ell),N^+_{\flow{s}w}(\ell'))&
  =\int_{-\infty}^0d(\ell(t),\ell'(t))\;e^{2\,t}\;dt+
  \int_0^{+\infty} d(\ell(0),\ell'(0))\;e^{-2\,t}\;dt\\ &
  \leq d(\ell,\ell')+\frac{1}{2}\,d(\ell(0),\ell'(0))\;.
\end{align*}
Since the footpoint map $\pi:\gengeod X\ra X$ defined by
$\ell\mapsto \ell(0)$ is $\frac{1}{2}$-Hölder-continuous (see
\cite[Prop.~3.2]{BroParPau19}), this proves that $N^+_{\flow{s}w}$ is
uniformly continuous (actually $\frac{1}{2}$-Hölder-continuous),
uniformly in $s$.

Conversely, note that by convexity, for all $\ell,\ell'\in
\wss(\flow{s}w)$, since $\ell_+=\ell'_+$, we have $d(\ell(t),\ell'(t))
\leq d(\ell(0),\ell'(0))$ for every $t\geq 0$. Hence
\begin{align*}
  d(\ell,\ell')&=\int_{-\infty}^0d(\ell(t),\ell'(t))\;e^{2\,t}\;dt+
  \int_0^{+\infty} d(\ell(t),\ell'(t))\;e^{-2\,t}\;dt\\&\leq
  \int_{-\infty}^0d(\ell(t),\ell'(t))\;e^{2\,t}\;dt+
  \int_0^{+\infty} d(\ell(0),\ell'(0))\;e^{-2\,t}\;dt=
  d(N^+_{\flow{s}w}(\ell),N^+_{\flow{s}w}(\ell'))\;.
\end{align*}
Therefore $(N^+_{\flow{s}w})^{-1}$ is $1$-Lipschitz, hence uniformly
continuous, uniformly in $s$.

\medskip
(2) Again since the footpoint map is
$\frac{1}{2}$-Hölder-continuous, there exists a constant $c>0$ such
that for every $\epsilon\in\;]0,1]$, for all $s,s' \in\RR$ and
$\ell,\ell'\in\wss(w)$, if $d(\flow s\ell, \flow{s'}\ell') \leq
\epsilon$, then $d(\ell(s),\ell'(s'))\leq c\,\epsilon^{\frac{1}{2}}$.
We may assume that $s\leq s'$.
\begin{center}
  \begin{picture}(0,0)%
\includegraphics{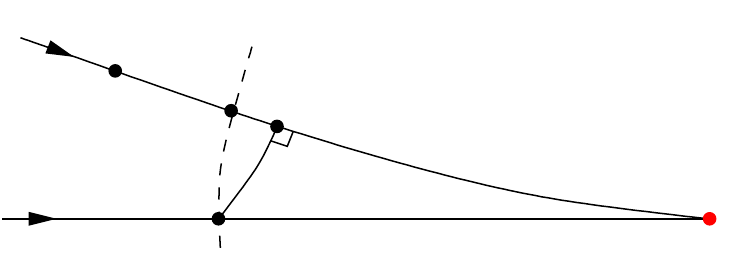}%
\end{picture}%
\setlength{\unitlength}{3812sp}%
\begingroup\makeatletter\ifx\SetFigFont\undefined%
\gdef\SetFigFont#1#2#3#4#5{%
  \reset@font\fontsize{#1}{#2pt}%
  \fontfamily{#3}\fontseries{#4}\fontshape{#5}%
  \selectfont}%
\fi\endgroup%
\begin{picture}(3618,1334)(439,-1220)
\put(1561,-1132){\makebox(0,0)[lb]{\smash{{\SetFigFont{11}{13.2}{\rmdefault}{\mddefault}{\updefault}{\color[rgb]{0,0,0}$\ell'(s')$}%
}}}}
\put(1861,-448){\makebox(0,0)[lb]{\smash{{\SetFigFont{11}{13.2}{\rmdefault}{\mddefault}{\updefault}{\color[rgb]{0,0,0}$p$}%
}}}}
\put(4042,-994){\makebox(0,0)[lb]{\smash{{\SetFigFont{11}{13.2}{\rmdefault}{\mddefault}{\updefault}{\color[rgb]{0,0,0}$\ell_+={}\ell'_+$}%
}}}}
\put(602, -9){\makebox(0,0)[lb]{\smash{{\SetFigFont{11}{13.2}{\rmdefault}{\mddefault}{\updefault}{\color[rgb]{0,0,0}$\ell$}%
}}}}
\put(998,-138){\makebox(0,0)[lb]{\smash{{\SetFigFont{11}{13.2}{\rmdefault}{\mddefault}{\updefault}{\color[rgb]{0,0,0}$\ell(s)$}%
}}}}
\put(508,-1165){\makebox(0,0)[lb]{\smash{{\SetFigFont{11}{13.2}{\rmdefault}{\mddefault}{\updefault}{\color[rgb]{0,0,0}$\ell'$}%
}}}}
\put(1202,-552){\makebox(0,0)[lb]{\smash{{\SetFigFont{11}{13.2}{\rmdefault}{\mddefault}{\updefault}{\color[rgb]{0,0,0}$\ell(s')$}%
}}}}
\end{picture}%

\end{center}
Since $\ell_+=\ell'_+$, by the convexity of the horoballs and by the
fact that closest point maps on nonempty closed convex subsets do not
increase the distances, with $p$ the closest point to $\ell'(s')$ on
$\ell([s,+\infty[)$, we have $p\in\ell([s',+\infty[)$ and
$$
|s-s'|=d(\ell(s),\ell(s'))\leq d(\ell(s),p)\leq
d(\ell(s),\ell'(s'))\leq c\,\epsilon^{\frac{1}{2}}\;.
$$
Let us fix $T>0$ and let us assume that $s\in[-T,T]$. By
\cite[Eq.~(2.8)]{BroParPau19}, we have $d(\flow{s'-s}\ell',\ell')
\leq |s-s'|$. By the change of variable $t\mapsto t+s$ in Equation
\eqref{eq:defdBartelLuck}, we have
$$
d(\ell,\flow{s'-s}\ell')\leq e^{2|s|}d(\flow s\ell,\flow{s'}\ell')\;.
$$
Therefore,
$$
d(\ell,\ell')\leq   d(\ell,\flow{s'-s}\ell')+  d(\flow{s'-s}\ell',\ell')
\leq e^{2T}\,\epsilon+c\,\epsilon^{\frac{1}{2}}\;.
$$

Conversely, for all $\epsilon\in\;]0,1]$, $T>0$, $s,s' \in[-T,T]$ and
$\ell,\ell'\in\wss(w)$, assume that $\max\{|s-s'|,\;d(\ell,\ell')\}\leq
\epsilon$.  Then by similar arguments, we have
$$
d(\flow s\ell,\flow{s'}\ell')\leq
d(\flow s\ell,\flow{s}\ell')+d(\flow s\ell',\flow s\flow{s'-s}\ell')
\leq e^{2T}(d(\ell,\ell')+|s'-s|)\leq 2\,e^{2T}\,\epsilon\;.
$$
This proves Assertion (2) of Lemma \ref{lem:unifcontwss}.

\medskip (3) Let $T>0$ and $s\in[-T,T]$. By Assertion (1), by the
commutativity of the diagram \eqref{eq:diagcomNplus} and by the
invertibility of $\flow{s}$, we only have to prove that $\flow{s}:\G
X\ra \G X$ is uniformly continuous, uniformly in $s\in[-T,T]$. As
already seen, for all $\ell,\ell'\in\G X$, we have $d(\flow
s\ell,\flow{s'}\ell')\leq e^{2T}\,d(\ell,\ell')$, hence the result
follows.
\cqfd

\bigskip 
We refer to \cite[\S7.2]{BroParPau19} for the following
definitions.  Let $\D^-=(D_i)_{i\in I^-}$ be a locally finite
$\Ga$-equivariant family of nonempty proper closed convex subsets of
$X$ and let $\D^+=(H_j)_{j\in I^+}$ be a locally finite
$\Ga$-equivariant family of (closed) horoballs in $X$. Let
$\sim=\sim_{+}$ be the equivalence relation on $I^+$ defined by $j\sim
j'$ if and only if $H_{j'}=H_j$ and there exists $\ga\in\Ga$ such that
$j'=\ga j$. Let $\sim=\sim_{-}$ be the similarly defined equivalence
relation on $I^-$.

For all $j\in I^+$ and $s\in \RR$, let $H_{j,s}$ be the horoball
contained in $H_j$ consisting of points at a distance at least $s$
from the complement of $H_j$ if $s\geq 0$, and otherwise, let
$H_{j,s}$ be the closed $(-s)$-neighbourhood $\N_{-s}H_{j}$ of $H_{j}$,
which is the horoball containing $H_j$ consisting of the points that
are at distance at most $-s$ from $H_j$.

\medskip
For every $j\in I^+$, let $w_j$ be any geodesic ray starting from the
boundary of the horoball $H_j$ and converging to the point at infinity
of $H_j$, so that $H\!B_+(w_j)=H_j$. We denote
\begin{align*}
&W^+_j=\wss(w_j),\;\; W^{0+}_j =\ws(w_j),\;\;
W^{0+}_{\sigma,j} =\ws_\sigma(w_j),\\& N^+_j=N^+_{w_j},
\;\;\mu^+_{j}=\muss{w_j},\;\;h_{j}=h_{w_j}
\;\;\;{\rm and}\;\;\;\mu^{0+}_{j}=\mus{w_j}\,.
\end{align*}
Using the homeomorphism $h_{j}$ from $W^+_j\times\RR$ to $W^{0+}_j$
defined by $(\ell,s)\mapsto \flow s \ell$ and the homeomorphism
$N^+_j:W^+_j\ra \normalin H_j$ defined by $\ell\mapsto \ell_{\mid\,
]-\infty,0]}$, for all $s\in\RR$ and $\ell\in W^+_j$, we thus have by
    Equations \eqref{eq:stabmeas} and \eqref{eq:strongstabmeas}
\begin{equation}\label{eq:strongstabmeasj}
  d\mu^{0+}_j(\flow s\ell)=
  e^{-\delta_{\Ga} s}\,d\,\wt\sigma^-_{H_j}(\ell_{\mid \,]-\infty,0]})\,ds\;.
\end{equation}

For all $j\in I^+$ and $s_0\in\RR$, since $H_j$ is the
$s_0$-neighbourhood of $H_{j,s_0}$ if $s_0\geq 0$ and since $H_{j,s_0}$
is the $(-s_0)$-neighbourhood of $H_{j}$ if $s_0\leq 0$, by
\cite[Eq.~(7.7)]{BroParPau19} (see also \cite[Prop.~4
  (iii)]{ParPau14ETDS} in the manifold case), for every $\ell\in
W^+_j$, we have
\begin{equation}\label{eq:skinningvoisinage}
  d\,\wt\sigma^-_{H_j}(\ell_{\mid \,]-\infty,0]})=e^{\delta_{\Ga} s_0}\,
  d\,\wt\sigma^-_{H_{j,s_0}}((\flow {s_0}\ell)_{\mid \,]-\infty,0]})\;.
\end{equation}

For every $t_0\in \RR$ fixed, we also define
\begin{equation}\label{eq:defimesuresomme}
\wt\sigma^{+}_{\D^-}=\sum_{i\in I^-/_\sim} \wt\sigma^{+}_{D_i}
\;\;\;{\rm and}\;\;\;\wt\mu^{0+}_{\D^+, t_0}=
\sum_{j\in I^+/_\sim} {\mu^{0+}_j}_{\mid W^{0+}_{t_0,j}}\;.
\end{equation}
Since the $\Ga$-equivariant family $(H_{j})_{j\in I^+}$ is locally
finite and since $t_0>-\infty$, the two measures $\wt\sigma^{+}_{\D^-}$
and $\wt\mu^{0+}_{\D^+, t_0}$ are locally finite. This is the reason
why it is important to restrict the (weak) stable leaves $W^{0+}_{j}$ to
their upper parts $W^{0+}_{t_0,j}$. These two measures are also
$\Ga$-equivariant by Equations \eqref{eq:equivskinn} and
\eqref{eq:equivweakstabmeas} (and by the $\Ga$-equivariance of the
families $\D^\pm$). Hence (see for instance \cite[\S 2.8]{PauPolSha15}
for the definition of the induced measure when $\Ga$ may have
torsion), they induce locally finite measures $\sigma^{+}_{\D^-}$ and
$\mu^{0+}_{\D^+, t_0}$ on $\Ga\bs\gengeod X$.


\section{Joint partial equidistribution of common perpendiculars
  to shrinking horoballs at a given density}
\label{sec:shrinking}

In this section, we prove, as an application of
\cite[Thm.~11.3]{BroParPau19}, a joint partial equidistribution
theorem for pairs consisting of a common perpendicular between a
locally convex subset and a quotient horoball on the one hand and its
image by the geodesic flow at a large time on the other hand. This
gives a generalised geometric version in negative curvature (including
variable one and in any dimension) of the case $n=2$ of
\cite[Thm.~6]{Marklof10Inv} and \cite[Thm.~6.1]{Lutsko22} for
hyperbolic surfaces.

With the notation of Section \ref{sec:background}, under the finiteness and
mixing property of the Bowen-Margulis measure and the finiteness and
nonvanishing property of the skinning measure, the image
$\flow{t}\Ga\normalout D_i$ by the geodesic flow at time $t$ of the
image in $\Ga\bs \G X$ of the outer normal bundle of $D_i$ (endowed
with its skinning measure) equidistributes towards the Bowen-Margulis
measure in $\Ga\bs\G X$. For a proof, we refer to  
\cite[Thm.~1]{ParPau14ETDS} in the manifold case and to \cite[Thm.~10.2
  with potential 0]{BroParPau19} in general.
We will take on $\flow{t}\Ga\normalout D_i$
sufficiently many images by $\flow{t}$ of common perpendiculars
from $D_i$ to $H_j$ in order to have a constant density with respect
to the skinning measure of $D_i$ in $\Ga\bs\G X$.

\medskip
For all $i\in I^-$ and $j\in I^+$ such that the point at infinity of
$H_{j}$ is not contained in $\partial_\infty D_i$ (or equivalently
such that $\partial_\infty D_i\,\cap \partial_\infty H_j=\emptyset$),
let $\rho_{i,j}$ be the unique geodesic ray in $\normalout D_i$ such
that $\rho_{i,j}(+\infty)$ is the point at infinity of $H_{j}$, and
let $\lambda_{i,j} = d(D_i,H_j)$.

\btheo \label{theo:Marklof6} Let $X$ be either a proper $\RR$-tree
without terminal points or a complete simply connected Riemannian
manifold with pinched negative sectional curvature at most $-1$. Let
$\Ga$ be a nonelementary discrete group of isometries of $X$.  Let
$\D^-=(D_i)_{i\in I^-}$ be a locally finite $\Ga$-equivariant family
of nonempty proper closed convex subsets of $X$ and let
$\D^+=(H_j)_{j\in I^+}$ be a locally finite $\Ga$-equivariant family
of horoballs in $X$.  Assume that the Bowen-Margulis measure $m_{\rm
  BM}$ on $\Ga\bs\G X$ is finite and mixing for the geodesic flow on
$\Ga\bs\G X$.  Then, for every $t_0\in\RR$, for the weak-star
convergence of measures on $\G_{+,0} X\times\gengeod X$, we have
\begin{equation}\label{eq:mainup}
\lim_{t\ra+\infty} \;\|m_{\rm BM}\|\;e^{-\delta_{\Ga}\, t}
\sum_{\substack{i\in I^-/_\sim,\; j\in I^+/_\sim, \;\ga\in\Ga\\ 
\partial_\infty D_i\,\cap \,\partial_\infty H_{\ga j}=\emptyset,\; 
\lambda_{i,\,\ga j}\leq t-t_0}} \; 
\Delta_{\rho_{i,\ga j}}\otimes\Delta_{\flow t\rho_{\ga ^{-1}i,j}}
\;=\wt\sigma^+_{\D^-}\otimes\wt \mu^{0+}_{\D^+,t_0}\;.
\end{equation}
\etheo

\dem Let us first give some notation that will be useful in this proof. For
all $s\in\RR$ and $(i,j)$ in $I^-\!\times I^+$ such that the closures
$\overline{D_i}$ and $\overline{H_{j,s}}$ of $D_i$ and $H_{j,s}$ in
$X\cup\partial_\infty X$ have empty intersection, let $\lambda_{i,j,s}
=d(D_i,H_{j,s})$ be the length of the common perpendicular from $D_i$
to $H_{j,s}$, and let $\alpha^-_{i,\,j,\,s} \in\gengeod X$ be its
parametrisation: it is the unique map from $\RR$ to $X$ such that

$\bullet$~ $\alpha^-_{i,\,j,\,s}(t)= \alpha^-_{i,\,j,\ s}(0)\in D_i$
if $t\leq 0$,

$\bullet$~ $\alpha^-_{i,\,j,\,s}(t)= \alpha^-_{i,\,j,\,s} (\lambda_{i,\,j,\,s})
\in H_{j,s}$ if $t\geq \lambda_{i,\,j,\,s}$, and 

$\bullet$~ ${\alpha^-_{i,j,s}}|_{ [0,\,\lambda_{i,j,s}]} =\alpha_{i,\,j,\,s}$
is the shortest geodesic arc starting from a point of $D_i$ and
ending at a point of $H_{j,s}$.  

\noindent
We have $\lambda_{i,j} = 0$ if $\overline{D_i}\,\cap\,
\overline{H_{j}}\neq\emptyset$ and $\lambda_{i,j} = \lambda_{i,j,0}>0$
if $\overline{D_i}\,\cap\, \overline{H_{j}}=\emptyset$, so that
$\lambda_{i,j,s}= \lambda_{i,j} +s$ when both terms $\lambda_{i,j}$
and $\lambda_{i,j,s}$ are defined and positive. Note that
$\lambda_{i,\ga j,s}=\lambda_{\ga^{-1}i,j,s}$ for every $\ga\in\Ga$,
by equivariance.  Let $\alpha^+_{i,j,s}= \flow{\lambda_{i,j,s}}
\alpha^-_{i,j,s}\in \gengeod X$, which is isometric exactly on
$[-\,\lambda_{i,j,s},0]$.

\medskip
The term on the left in Equation \eqref{eq:mainup} is independent of
the choice of the representatives of $i$ and $j$. Let us fix $(i,j)\in
I^-\times I^+$ and let us prove that for the weak-star convergence of
measures on $\G_{+,0} X\times\gengeod X$, we have
\begin{equation}\label{eq:formuleamontrer}
\lim_{t\ra+\infty} \;\|m_{\rm BM}\|\;e^{-\delta_{\Ga}\, t}\!\!\!
\sum_{\substack{\ga\in\Ga
    \\\partial_\infty D_i\,\cap\, \partial_\infty H_{\ga j}=\emptyset,
\;\lambda_{i,\,\ga j}\leq t-t_0}} \; 
\Delta_{\rho_{i,\ga j}}\otimes\Delta_{\flow t\rho_{\ga ^{-1}i,j}}
\;=\wt\sigma^+_{D_i}\otimes{\mu^{0+}_{j}}_{\mid W^{0+}_{t_0,j}}\;.
\end{equation}
The result follows by a (locally finite) summation using Equation
\eqref{eq:defimesuresomme}.

\medskip
For all $\tau\in\;]0,1]$ and $s_0\ge t_0$, Theorem 11.3 of
\cite{BroParPau19} (in the case with  potential $0$) applied to
the locally finite $\Ga$-equivariant families $(D_{\alpha
  i})_{\alpha\in\Ga}$ and $(H_{\beta j,s_0})_{\beta\in\Ga}$ (see also
\cite[Eq.~(12)]{ParPau17ETDS} in the manifold case) gives, for the
weak-star convergence of measures on $\gengeod X\times\gengeod X$,
\begin{equation}\label{eq:formuledontonpart}
\lim_{t\ra+\infty} \;
\|m_{\rm BM}\|\;e^{-\delta_{\Ga}\, t}\sum_{\substack{\ga\in\Ga,\; 
\ov{D_i}\,\cap \,\ov{H_{\ga j,s_0}}\,=\emptyset\\ 
t-\tau<\lambda_{i,\,\ga j,\,s_0}\leq t}} \; 
\Dirac_{\alpha^-_{i,\,\ga j,\,s_0}} \otimes\Dirac_{\alpha^+_{\ga^{-1}i,\,j,\,s_0}}
=\frac{1-e^{-\delta_{\Ga}\,\tau}}{\delta_{\Ga}} \,
\wt\sigma^+_{D_i}\otimes \wt\sigma^-_{H_{j,s_0}}\;.
\end{equation}

Let us consider two Borel subsets $\Omega^-$ of $\normalout D_i$ and
$\Omega^+$ of $W^+_j$ with positive and finite measure for
$\wt\sigma^+_{D_i}$ and $\mu^{+}_{j}$ respectively, whose boundaries
have zero measure for $\wt\sigma^+_{D_i}$ and $\mu^{+}_{j}$
respectively. For all $s_0\geq t_0$ and $\tau>0$, the product
$B=\Omega^-\times h_j(\Omega^+\times[s_0,s_0+\tau])$ is contained in
$\normalout D_i\times W^{0+}_{t_0,j}$.

\medskip
\noindent {\bf Step 1. } Let us first relate the two right hand sides
of Equations \eqref{eq:formuleamontrer} and
\eqref{eq:formuledontonpart} evaluated on the Borel set $B$.

By respectively Equation \eqref{eq:strongstabmeasj}, Equation
\eqref{eq:skinningvoisinage}, an easy integral computation and the
commutativity of the diagram \eqref{eq:diagcomNplus}, we have
\begin{align}
&\int_{(\rho,\ell,s)\in \Omega^-\times\Omega^+\times[s_0,s_0+\tau]}
d\,\wt\sigma^+_{D_i}(\rho)\,d\mu^{0+}_{j}(\flow s\ell)\nonumber\\=\;&
\int_{(\rho,\ell,s)\in \Omega^-\times\Omega^+\times[s_0,s_0+\tau]}
d\,\wt\sigma^+_{D_i}(\rho)\,e^{-\delta_{\Ga} s}d\wt\sigma^-_{H_j}
(\ell_{\mid \,]-\infty,0]})\,ds\nonumber\\=\;&
  \int_{(\rho,\ell)\in \Omega^-\times\Omega^+}
  d\,\wt\sigma^+_{D_i}(\rho)\,\Big(\int_{s_0}^{s_0+\tau}
  e^{-\delta_{\Ga} s}e^{\delta_{\Ga} s_0}\, ds\Big)
  d\,\wt\sigma^-_{H_{j,s_0}}((\flow {s_0}\ell)_{\mid \,]-\infty,0]})
      \nonumber\\=\;&\int_{(\rho,\ell)\in \Omega^-\times\Omega^+}
  \frac{1-e^{-\delta_{\Ga}\,\tau}}{\delta_{\Ga}}\;d\wt\sigma^+_{D_i}(\rho)\,
  d\,\wt\sigma^-_{H_{j,s_0}}((\flow {s_0}\ell)_{\mid \,]-\infty,0]})
      \nonumber\\=\;&
  \int_{(\rho,\rho')\in \Omega^-\times\flow{s_0}_\mid N^+_j(\Omega^+)}
  \frac{1-e^{-\delta_{\Ga}\,\tau}}{\delta_{\Ga}}\;d\,\wt\sigma^+_{D_i}(\rho)\,
  d\,\wt\sigma^-_{H_{j,s_0}}(\rho')\;.\label{eq:travailmebredroite}
\end{align}

\medskip
\noindent {\bf Step 2. } Let us now relate the two index sets of the
left hand sides of Equations \eqref{eq:formuleamontrer} and
\eqref{eq:formuledontonpart}, except for the ranges of
$\lambda_{i,\,\ga j}$ and $\lambda_{i,\,\ga j,\,s_0}$, that will be
taken care of in Step 3.

For every $\ga\in\Ga$, if $\ov{D_i}\,\cap \,\ov{H_{\ga j,s_0}} \,=
\emptyset$ (so that $\alpha^-_{i,\,\ga j,\,s_0}$ is defined), then
$\partial_\infty D_i\cap \partial_\infty H_{\ga j} =\emptyset$ (so
that $\rho_{i,\ga j}$ is defined) and $\alpha^-_{i,\,\ga
  j}(0)=\rho_{i,\,\ga j}(0)$. Conversely, we have the following
result, whose Case 2 will only by useful for the proof of the narrow
convergence claim in Theorem \ref{theo:quotient}. Recall that
$\pi:\gengeod X\ra X$ is the footpoint projection.

\blemm\label{lem:dichot} If $\partial_\infty D_i\cap \partial_\infty
H_{\ga j} =\emptyset$ (so that $\rho_{i,\ga j}$ is defined) and if

{\bf Case 1. } either the Borel sets $\Omega^\pm$ are relatively compact,

{\bf Case 2. } or the stabilizer $\Ga_{D_i}$ of $D_i$ in $\Ga$ acts
cocompactly on $\partial D_i$,

\noindent
then there exists a finite subset $F$ of $\Ga$ (depending on
$i,j,\Omega_-,t_0$), such that for all $\ga\in \Ga-F$ and $s_0\geq
t_0$, if $\rho_{i,\ga j}(0)\in \pi(\Omega_-)$, then $\ov{D_i}\,\cap
\,\ov{H_{\ga j,s_0}} \, = \emptyset$ (so that $\alpha^-_{i,\,\ga
  j,\,s_0}$ is defined).
\elemm

\dem This follows directly by the local finiteness of the family
$(H_j)_{j\in I_+}$ in Case 1, and by this local finiteness and the
fact that the skinning measure on $\partial^1_+D_i$ is invariant under 
$\Ga_{D_i}$ in Case 2. \cqfd

\medskip
\noindent {\bf Step 3. } Let us finally relate the two pairs of Dirac
masses in the left hand sides of Equations \eqref{eq:formuleamontrer}
and \eqref{eq:formuledontonpart}, as well as the ranges of
$\lambda_{i,\,\ga j}$ and $\lambda_{i,\,\ga j,\,s_0}$.

If $\ga\in\Ga-F$ furthermore satisfies $\lambda_{i,\,\ga j}\geq T$ for
some $T>0$ (which excludes only finitely many more $\ga\in\Gamma$),
then the generalised geodesics $\rho_{i,\ga j}$ and $\alpha^-_{i,\,\ga
  j,\,s_0} $ coincide on $]-\infty, T+s_0]$, hence on $]-\infty,
T+t_0]$. Therefore, they are at distance at most $\epsilon$ for any
given $\epsilon>0$ if $T$ is large enough (uniformly in $s_0$ and
$\ga$) by Equation \eqref{eq:defdBartelLuck}.

Since $X$ has extendible geodesics, for every $\ga\in\Ga$ such that
$\partial_\infty D_i\cap \partial_\infty H_{\ga j} =\emptyset$ (or
equivalently $\partial_\infty D_{\ga^{-1}i}\cap \partial_\infty H_{j}
=\emptyset$), let $\wt\rho_{\ga^{-1}i,j}\in\G X$ be any geodesic line 
such that we have $\wt\rho_{\ga^{-1}i,j}\!\!  \mid_{[0,+\infty[} \,=
\rho_{\ga^{-1}i,j} \!\!  \mid_{[0,+\infty[}$. For $t$ large enough,
the generalised geodesics $\flow t\wt\rho_{\ga^{-1}i,j}$ and $\flow t 
\rho_{\ga^{-1}i,j}$, which coincide on $[-t,+\infty[$, are arbitrarily
close (uniformly in $\ga$) by Equation \eqref{eq:defdBartelLuck}.
Hence we may replace $\flow t \rho_{\ga^{-1}i,j}$ by $\flow t
\wt\rho_{\ga^{-1}i,j}$ in the formula \eqref{eq:formuleamontrer} that
we want to prove.

Note that $\flow t\wt\rho_{\ga^{-1}i,j}$ belongs to $W^{0+}_j$, and
that $\flow {\lambda_{i,\ga j}}\wt\rho_{\ga^{-1}i,j}$ belongs to
$W^{+}_j$. Since
$$
\flow t\wt\rho_{\ga^{-1}i,j}= \flow {t-\lambda_{i,\ga j }}
(\flow {\lambda_{i,\ga j }}\wt\rho_{\ga^{-1}i,j })\;,
$$
it follows from Lemma \ref{lem:unifcontwss} (2) that $\flow
t\wt\rho_{\ga^{-1}i,j}$ is close to the subset $h_j(\Omega^+\times
[s_0,s_0+\tau])$ if and only if $t-\lambda_{i,\ga j }$ is close to
$[s_0,s_0+\tau]$ and $\flow {\lambda_{i,\ga j}}\wt\rho_{\ga^{-1}i,j}$
is close to $\Omega^+$.  In particular, if $\flow t \wt
\rho_{\ga^{-1}i,j}$ is close to $h_j(\Omega^+\times [s_0,s_0+\tau])$
and $t$ is large enough, then $\lambda_{i,\ga j }$ is large enough,
and $\lambda_{i,\ga j, s_0}$ is close to $[t-\tau,t]$ (uniformly in
$\ga$).

Finally, the negative geodesic ray $\flow{s_0}_\mid N^+_j(\flow
{\lambda_{i,\ga j}}\wt\rho_{\ga^{-1}i,j})$, which is close to the
subset $\flow{s_0}_\mid N^+_j(\Omega^+)$ by Lemma
\ref{lem:unifcontwss} (1) and (3), coincides with the generalized
geodesic $\alpha^+_{\ga^{-1}i,j,s_0}$ on the whole interval
$]-\lambda_{i,\ga j,s_0}, +\infty[$. Since $\lambda_{i,\ga j, s_0}$ is
    large (uniformly in $\ga$) when $t$ is large, and again by
    Equation \eqref{eq:defdBartelLuck}, this implies that the
    generalised geodesic lines $\flow{s_0}_\mid N^+_j(\flow
    {\lambda_{i,\ga j }}\wt\rho_{\ga^{-1}i,j})$ and
    $\alpha^+_{\ga^{-1}i,j,s_0}$ are close (uniformly in $\ga$).

\medskip
In order to conclude the proof of the weak-star convergence in Theorem
\ref{theo:Marklof6}, we assume that $\Omega^\pm$ are relatively
compact (so that Case 1 of Lemma \ref{lem:dichot} applies). We
evaluate the two sides of Formula \eqref{eq:formuledontonpart} on the
relatively compact Borel subset $\Omega^-\times\flow{s_0}_\mid
N^+_j(\Omega^+)$ of $\gengeod X\times\gengeod X$, whose boundary has
measure zero for the limit measure. By using Formula
\eqref{eq:travailmebredroite}, this implies Formula
\eqref{eq:formuleamontrer} evaluated on the relatively compact Borel
subset $\Omega^-\times\Omega^+\times[s_0,s_0+\tau]$, whose boundary
has measure zero for the limit measure. The result follows.  \cqfd

\bigskip
Let us now give a version of Theorem \ref{theo:Marklof6} in the
discrete tree case.  Let $\XX$ be a locally finite simplicial tree
without terminal vertices, with geometric realisation $X=|\XX|_1$
(with edge lengths equal to $1$) and with boundary at infinity
$\partial_\infty\XX=\partial_\infty X$. We denote by $V\XX$ the set of
vertices of $\XX$, identified with its image in $X$. Let $\Ga$ be a
nonelementary discrete subgroup of the inversion-free automorphism
group $\Aut(\XX)$ of $\XX$, and let $\delta_\Ga>0$ be its critical
exponent. We refer to \cite[\S 2.6]{BroParPau19} for background, as
well as for the definition of the space of generalised discrete
geodesic lines
\[
\gengeod\XX=\{\ell\in\gengeod X: \ell(0)\in
V\XX,\;\ell_\pm\in V\XX\cup\partial_\infty \XX\}
\]
of $\XX$, and the definition of the discrete-time geodesic flow
$(\flow n)_{n\in\ZZ}$ on $\gengeod\XX$, given by setting $\flow
n\ell:m\mapsto \ell(m+n)$ for all $\ell\in\gengeod \XX$ and
$m,n\in\ZZ$.

By taking the intersections with $\gengeod \XX$ of the previously
defined objects for $X$, we define (see op.~cit.)~the closed subspaces
$\G\XX$, $\G_\pm\XX$ and $\G_{\pm,0}\XX$ of $\gengeod \XX$, the strong
stable leaf $W^+(w)$, the stable leaf $W^{0+}(w)$ and the truncated
stable leaf
\[
W^{0+}_{n_0}(w)= \bigcup_{n\in\ZZ, n\geq n_0}\flow n W^+(w)
\]
of $w\in\G_+\XX$, where $n_0\in\ZZ$, and the outer/inner unit normal
bundles $\normalpm\DD$ of a nonempty proper simplicial subtree
$\DD$ of $\XX$. We define similarly (see op.~cit.) the outer/inner
skinning measure $\wt \sigma^\pm_\DD$ on $\normalpm\DD$ and the
Bowen-Margulis measures $\wt m_{\rm BM}$ on $\G\XX$ and $m_{\rm BM}$
on $\Ga\bs\G\XX$ associated with any choice of Patterson-Sullivan
density $(\mu_x)_{x\in V\XX}$.

Given $w\in\G_+\XX$, its stable horoball $H\!B_+(w)$ is a subtree of
$\XX$ and we again denote by $N^+_w: \wss(w)\to\normalin H\!B_+(w)$
the canonical homeomorphism defined in Section \ref{sec:background}.
We now have a homeomorphism $h_w:\wss(w)\times\ZZ\ra \ws(w)$ defined
by $(\ell,m)\mapsto \flow m\ell$. The conditional measure $\mus{w}$ of
the Bowen-Margulis measure $\wt m_{\rm BM}$ (for the discrete-time
geodesic flow on $\gengeod \XX$) on the stable leaf $\ws(w)$ of $w$ is
now defined, for $m\in\ZZ$ and $\ell\in \wss(w)$, by
\begin{equation}\label{eq:stabmeasdiscr}
d\mus{w}(\flow m\ell)=e^{-\delta_\Ga m}d\muss{w}(\ell)\,dm\;,
\end{equation}
where $dm$ is the counting measure on $\ZZ$.

\medskip 
Let $\D^-= (\DD^-_i)_{i\in I^-}$ and $\D^+=(\HH^+_j)_{j\in I^+}$ be
locally finite $\Ga$-equivariant families of nonempty proper
simplicial subtrees of $\XX$, with $\HH^+_j$ a horoball for every
$j\in I^+$. We consider the geometric realisations $D_i=|\DD_i|_1$ of
$\DD_i$ and $H_j=|\HH_j|_1$ of $\HH_j$. For every $n_0\in\ZZ$, we
define the horoball $H_{j,n_0}$ such that $H_j$ is the
$n_0$-neighbourhood of $H_{j,n_0}$ if $n_0\geq 0$ and $H_{j,n_0}$ is
the $(-n_0)$-neighbourhood of $H_{j}$ if $n_0\leq 0$. For every
$n_0\in\ZZ$, as in the end of Section \ref{sec:background}, we define
the measures $\wt\sigma^{+}_{\D^-}$ and $\wt\mu^{0+}_{\D^+, n_0}$ on
$\gengeod \XX$, and their induced measures $\sigma^{+}_{\D^-}$ and
$\mu^{0+}_{\D^+, n_0}$ on $\Ga\bs \gengeod \XX$.

For all $m\in\ZZ$ and $(i,j)\in I^-\times I^+$, the elements
$\rho_{i,\,j}$ and $\alpha^\pm_{i,\,j,\,m}$, respectively defined just
before and just after the statement of Theorem \ref{theo:Marklof6},
actually belong to $\gengeod \XX$.

Note that for many interesting lattices in $\Aut(\XX)$ (and this will
turn out to be the case for the application in Subsection
\ref{subsec:nonarchFarey}), the time-one geodesic flow is not mixing
(it is not even ergodic), though the time-two geodesic flow is mixing
on an halfsubspace, see \cite[End of \S 4.4]{BroParPau19} for
explanations. This explains the usefulness of Assertion (2) in the
next statement.

Fix a basepoint $x^\bullet\in V\XX$.  Let $V_{\rm even}\XX$ be the
subset of $V\XX$ of vertices at even distance from $x^\bullet$. Let
$$
\gengeod_{\rm even} \XX=\{\ell\in\gengeod\XX : \ell(0)\in V_{\rm
  even}\XX\}\hspace{3mm}\textrm{and}\hspace{3mm} \G_{\rm even}\XX
=\gengeod_{\rm even}\XX\cap\G\XX\,.
$$

\btheo \label{theo:Marklof6discrete} Let $\XX$ be a locally finite
simplicial tree without terminal vertices. Let $\Ga$ be a
nonelementary discrete subgroup of $\Aut(\XX)$.  Let
$\D^-=(\DD^-_i)_{i\in I^-}$ and $\D^+=(\HH^+_j)_{j\in I^+}$ be locally
finite $\Ga$-equivariant families of nonempty proper simplicial
subtrees of $X$, with $\HH^+_j$ an horoball for every $j\in I^+$.

(1) Assume that the Bowen-Margulis measure $m_{\rm BM}$ on
$\Ga\bs\G\XX$ endowed with the discrete-time geodesic flow is
finite and mixing. Then, for every $n_0\in\ZZ$, for the weak-star
convergence of measures on $\G_{+,0} \XX\times\gengeod \XX$, we have
$$
\lim_{n\ra+\infty} \;\|m_{\rm BM}\|\;e^{-\delta_\Ga\, n}
\sum_{\substack{i\in I^-/_\sim,\; j\in I^+/_\sim, \;\ga\in\Ga\\ 
\partial_\infty \DD_i\,\cap \,\partial_\infty \HH_{\ga j}=\emptyset,\; 
\lambda_{i,\,\ga j}\leq n-n_0}} \; 
\Delta_{\rho_{i,\ga j}}\otimes\Delta_{\flow n\rho_{\ga ^{-1}i,j}}
\;=\wt\sigma^+_{\D^-}\otimes\wt \mu^{0+}_{\D^+,n_0}\;.
$$

(2) Assume that $\Ga$ preserves the $V_{\rm even}\XX$. Assume that the restriction to $\Ga\bs \G_{\rm even} \XX$ of
the Bowen-Margulis measure $m_{\rm BM}$ is finite and mixing for the
time-two of the discrete-time geodesic flow. Assume that the endpoints
of every common perpendicular between disjoint elements of $\D^-$ and
$\D^+$ belong to $V_{\rm even}\XX$.  Then, for every $n_0\in\ZZ$, for
the weak-star convergence of measures on $\G_{+,0} \XX\times\gengeod
\XX$, we have
\begin{align*}
\lim_{n\ra+\infty} &\;\frac{\|m_{\rm BM}\|}{2}\;e^{- 2\,\delta_\Ga\,n}
\sum_{\substack{i\in I^-/_\sim,\; j\in I^+/_\sim, \;\ga\in\Ga\\ 
\partial_\infty \DD_i\,\cap \,\partial_\infty \HH_{\ga j}=\emptyset,\; 
\lambda_{i,\,\ga j}\leq 2n-2n_0}} \; 
\Delta_{\rho_{i,\ga j}}\otimes\Delta_{\flow {2n}\rho_{\ga ^{-1}i,j}}
\\&=\wt\sigma^+_{\D^-}\,_{\mid \gengeod_{\rm even}\XX}
\otimes\,\wt \mu^{0+}_{\D^+,2n_0}\,_{\mid \gengeod_{\rm even}\XX}\;.
\end{align*}
\etheo

\dem (1) Let us fix $i\in I^-$ and $j\in I^+$. It follows from (the
case with zero potential of) \cite[Thm.~11.9]{BroParPau19} in the same
way as Thm.~11.3 of op.~cit.~follows from Thm.~11.1 ibid.~that for
every integer $m_0\geq n_0$, we have
\[
\lim_{n\ra+\infty} \;\|m_{\rm BM}\|\;e^{-\delta_\Ga\, n}
\sum_{\substack{\ga\in\Ga\\ {\DD^-_i}\cap {\HH^+_{\ga j,m_0}}
    =\emptyset,\; \lambda_{i,\,\ga j,m_0}= n}} \;
\Dirac_{\alpha^-_{i,\,\ga j}} \otimes\Dirac_{\alpha^+_{\ga^{-1}i,\,j,m_0}}
\;=\; \wt\sigma^+_{D_i}\otimes \wt\sigma^-_{H_{j,m_0}}\;
\]
for the weak-star convergence of measures on the locally compact space
$\gengeod \XX\times \gengeod \XX$. The proof of Theorem
\ref{theo:Marklof6discrete} (1) is then similar to that of Theorem
\ref{theo:Marklof6} using this equation instead of Equation
\eqref{eq:formuledontonpart}.

\medskip
(2) Let us fix $i\in I^-$ and $j\in I^+$. It follows from (the case
with zero potential of) now \cite[Thm.~11.11]{BroParPau19} (and more
precisely of Equation (11.28) in its proof) in the same way as
Thm.~11.3 of op.~cit.~follows from Thm.~11.1 ibid.~that for every
integer $m_0\geq n_0$, we have
\[
\lim_{n\ra+\infty} \;\frac{\|m_{\rm BM}\|}{2}\;e^{-2\,\delta_\Ga\, n}
\sum_{\substack{\ga\in\Ga\\{\DD^-_i}\cap {\HH^+_{\ga j,2m_0}}=
    \emptyset\\ \lambda_{i,\,\ga j,2m_0}= 2n}} \;
\Dirac_{\alpha^-_{i,\,\ga j}} \otimes\Dirac_{\alpha^+_{\ga^{-1}i,\,j,2m_0}}
\;=\; \wt\sigma^+_{D_i}\,_{\mid \gengeod_{\rm even}\XX}
\otimes \,\wt\sigma^-_{H_{j,2m_0}}\,_{\mid \gengeod_{\rm even}\XX}
\]
for the weak-star convergence of measures on the locally compact space
$\gengeod_{\rm even} \XX\times \gengeod_{\rm even} \XX$. The proof of
Theorem \ref{theo:Marklof6discrete} (2) is then similar to that of
Theorem \ref{theo:Marklof6} using this equation instead of Equation
\eqref{eq:formuledontonpart}.
\cqfd

\medskip
In order to conclude Section \ref{sec:shrinking}, let us give
equidistribution statements in the quotient by $\Ga$ of the two
previous results. In order to simplify them, we assume that $D$ is a
proper nonempty closed convex subset of $X$ and that $H$ is a (closed)
horoball of $X$ such that the $\Ga$-equivariant families $\D^-=(\ga
D)_{\ga\in\Ga}$ and $\D^+=(\ga H)_{\ga\in\Ga}$ are locally finite. In
the simplicial tree case as above, we assume that $D$ and $H$ are the
geometric realisations of simplicial subtrees $\DD$ and $\HH$ of
$\XX$.

We denote by $\Ga_D$ and $\Ga_H$ the stabilisers of $D$ and $H$ in
$\Ga$ respectively. For every $\ga\in\Ga$ such that the point at
infinity of $\ga H$ does not belong to $\partial_\infty D$, we define
the {\it multiplicity} of the common perpendicular from $D$ to $\ga H$
by
$$
m_\ga=\frac{1}{\card\big(\Ga_D\cap(\ga\Ga_H\ga^{-1})\big)}
$$
and we denote by $\rho_\ga$ the unique geodesic ray in $\normalout D$
converging to the point at infinity of $\ga H$. Note that for all
$\alpha\in\Ga_D$ and $\beta\in\Ga_H$, we have
$$
m_\ga=m_{\alpha\ga\beta}\quad{\rm and}\quad \alpha\rho_\ga=
\rho_{\alpha\ga\beta}\;.
$$

\btheo\label{theo:quotient} (1) For every $t_0\in\RR$, if
$(X,\Ga)$ satisfies the assumptions of Theorem \ref{theo:Marklof6} for
$\D^\pm$ as above, if furthermore the measures $\sigma^{+}_{\D^-}$ and
$\mu^{0+}_{\D^+, t_0}$ on $\Ga\bs\gengeod X$ are finite and nonzero,
then for the weak-star convergence of measures on $(\Ga\bs\G_{+,0}
X)\times(\Ga\bs\gengeod X)$, we have
\begin{equation}\label{eq:downmanifcase}
\lim_{t\ra+\infty} \;\|m_{\rm BM}\|\;e^{-\delta_\Ga\, t}
\sum_{\substack{\ga\in\Ga_{D}\!\backslash \Ga/\Ga_{H}\\ 
0<d(D,\ga H)\leq t-t_0}} \; 
m_\ga \;\Delta_{\Ga\rho_\ga}\otimes\Delta_{\flow t\Ga\rho_\ga}
\;=\sigma^+_{\D^-}\otimes\mu^{0+}_{\D^+,t_0}\;.
\end{equation}

(2) For every $n_0\in\ZZ$, if $(\XX,\Ga)$ satisfies the assumptions of
Theorem \ref{theo:Marklof6discrete} (1) for $\D^\pm$ as above, if
furthermore the measures $\sigma^{+}_{\D^-}$ and $\mu^{0+}_{\D^+,
  n_0}$ on $\Ga\bs\gengeod \XX$ are finite and nonzero, then for the
weak-star convergence of measures on $(\Ga\bs\G_{+,0}\XX)\times
(\Ga\bs\gengeod \XX)$, we have
$$
\lim_{n\ra+\infty} \;\|m_{\rm BM}\|\;e^{-\delta_\Ga\, n}
\sum_{\substack{\ga\in\Ga_{D}\!\backslash \Ga/\Ga_{H}\\ 
\partial_\infty D\,\cap \,\ga\partial_\infty H=\emptyset,\; 
d(D,\ga H)\leq n-n_0}} \; 
m_\ga \;\Delta_{\Ga\rho_\ga}\otimes\Delta_{\flow n\Ga\rho_\ga}
\;=\sigma^+_{\D^-}\otimes\mu^{0+}_{\D^+,n_0}\;.
$$

(3) For every $n_0\in\ZZ$, if $(\XX,\Ga)$ satisfies the assumptions of
Theorem \ref{theo:Marklof6discrete} (2) for $\D^\pm$ as above, if
furthermore the measures $\sigma^{+}_{\D^-}$ and $\mu^{0+}_{\D^+,
  n_0}$ on $\Ga\bs\gengeod \XX$ are finite and nonzero, then for the
weak-star convergence of measures on $(\Ga\bs\G_{+,0}\XX)\times
(\Ga\bs\gengeod \XX)$, we have
\begin{align*}
\lim_{n\ra+\infty} &\;\frac{\|m_{\rm BM}\|}{2}\;e^{-2\,\delta_\Ga\, n}
\sum_{\substack{\ga\in\Ga_{D}\!\backslash \Ga/\Ga_{H}\\ 
\partial_\infty D\,\cap \,\ga\partial_\infty H=\emptyset,\; 
d(D,\ga H)\leq 2n-2n_0}} \; 
m_\ga \;\Delta_{\Ga\rho_\ga}\otimes\Delta_{\flow{2n}\Ga\rho_\ga}
\\&=\sigma^+_{\D^-}\,_{\mid \;\Ga\bs\gengeod_{\rm even}\XX}\otimes
\mu^{0+}_{\D^+,2n_0}\,_{\mid \;\Ga\bs\gengeod_{\rm even}\XX}\;.
\end{align*}

Furthermore, when the stabilizer $\Ga_D$ of $D$ in $\Ga$ acts
cocompactly on $\partial D$, then the weak-star convergence claimed
in the three previous assertions may be improved to a narrow
convergence.
\etheo

\dem The first assertion follows from Theorem \ref{theo:Marklof6} in
the same way as Corollary 12.3 in the manifold case and Theorem 12.8
in the tree case of \cite{BroParPau19} follows from Theorem 11.1 of
\cite{BroParPau19}. The second and third assertions follow
respectively from Theorem \ref{theo:Marklof6discrete} (1) and (2) in
the same way as Theorems 12.9 and 12.12 of \cite{BroParPau19} follow
from Theorems 11.9 and 11.11 of \cite{BroParPau19}. The final narrow
convergence claim uses Case 2 of Lemma \ref{lem:dichot} and the fact
that $\Omega^\pm$ where only assumed to have nonzero and finite
skinning measures in the proof of Theorem \ref{theo:Marklof6}.
\cqfd

\medskip
\rem Assume first in this remark that $X$ is a (negatively curved)
symmetric space, that $\Ga$ is an arithmetic lattice and that $D$ has
smooth boundary. Note that the Bowen-Margulis measure is then the
Liouville measure, and in particular is a smooth measure. For all
$\ell\in\NN$ and $f\in\C^\ell_c(\Ga\bs T^1\wt M)$, we denote by
$\|f\|_\ell$ the Sobolev norm of $f$. We identify $\G_{+,0} X$ and
$\gengeod X$ with $T^1X$ by uniquely extending geodesic rays and
segments to geodesic lines. Then one could prove, as in \cite[Thm.~15
  (2)]{ParPau17ETDS} (see also \cite[Thm.~12.7 (2)]{BroParPau19}),
by replacing the above Equation \eqref{eq:formuledontonpart}
by the difference of the evaluations at $T=t$ and $T=t-\tau$ of
Equation (28) of \cite{ParPau17ETDS}, that there exists $\tau'>0$
such that we have an error term of the form $\bigO_{t_0}(e^{-\kappa'
  t}\|\psi^-\|_\ell \|\psi^+\|_\ell)$ when evaluating (before taking
the limit on the left hand side) the two sides of Equation
\eqref{eq:downmanifcase} on a pair of functions $\psi^\pm\in
\C^\ell_c(\Ga\bs T^1\wt M)$.

Assume now, with the notation of Section \ref{subsec:nonarchFarey},
that $X$ is the geometric realisation of the Bruhat-Tits tree $\XX_v$
of a $(\PGL_2,K_v)$ and $\Ga=\PGL_2(R_v)$ is the Nagao lattice. One
could prove a similar error term replacing a Sobolev regularity by a
locally constant regularity, as in Remark (ii) in \cite[page
  282]{BroParPau19} using \cite[Proposition 15.7]{BroParPau19} in
order to check the main assumption of that remark.

\section{Applications to equidistribution of Farey fractions}
\label{sec:applicationexamples}

In this section, we give five examples of applications of the results
of Section \ref{sec:shrinking}, by taking arithmetic families of
points (of Farey fractions type) with a given average density in an
expanding closed horophere, and we study their equidistribution
properties. As their proofs, though having similar schemes, make
reference to many different papers, and require numerous different
computations and checkings, it has not been possible, if only for the
sake of the readability of this paper, to regroup them into one
statement. More corollaries of Theorem \ref{theo:quotient} (1) may be
obtained by varying a nonuniform arithmetic lattice $\Ga$ in the
isometry group of a negatively curved symmetric space $X$. In Sections
\ref{subsec:retrouvonsMarklof} and \ref{subsec:complexFarey}, we
denote by $\begin{bmatrix} a & b \\ c & d \end{bmatrix}$ the image in
$\PSLC=\SLC/\{\pm\id\}$ of $\begin{pmatrix} a & b \\ c &
  d \end{pmatrix}\in \SLC$.

\subsection{Standard Farey fractions and Marklof's theorem}
\label{subsec:retrouvonsMarklof}

Let us now check that as a corollary of Theorem \ref{theo:quotient}
(1), we obtain a new and geometric proof of the case $n=2$ of
\cite[Thm.~6]{Marklof10Inv}.  We give extra details in the proof of
Corollary \ref{coro:marklof}, as it will serve as a model for the four
next examples.

Let $G=\PSLR$ and let $\Ga$ be the modular group $\PSLZ$. For all
$r,t\in \RR$,  let
\[
n_-(r)=\begin{bmatrix} 1 & r\\0&1\end{bmatrix}\quad 
{\rm and}
\quad 
\Phi^t=\begin{bmatrix} \;e^{-t/2} & 0\\0&e^{t/2}\end{bmatrix}\,.
\]
Let $$H=\{n_-(r):r\in\RR\}\,,$$ and let $$\Ga_H=H\cap \Ga=\{n_-(r):
r\in\ZZ\}\,.
$$
We see $\Ga_H\bs H$ as contained in $\Ga\bs G$, and we
endow $\Ga_H\bs H$ with its $H$-invariant probability measure
$\mu_{\Ga_H\bs H}$. We endow $\RR/\ZZ$ with its probability Haar
measure $dx$, so that the map $r\mapsto n_-(r)$ induces a measure
preserving homeomorphism $\RR/\ZZ\ra\Ga_H\bs H$.

For every $t\in\RR$, we consider the subset $\F_{t}$ of $\RR/\ZZ$
consisting of the {\em (standard) Farey fractions of height at most}
$e^{t/2}$, defined by
\[
\F_{t}=\Big\{\frac pq\!\!\mod 1:
\;p,q\in\ZZ,\quad (p,q)=1,\quad 0<q\leq e^{t/2}\Big\}\,.
\]
Note that both in the definition of $\Phi^t$ and $\F_{t}$, Marklof
replaces $t$ by $2t$, but our convention is more natural considering
the left part of Equation \eqref{eq:conjugflogeod} below.

Let $\Theta:\Ga\bs G\ra \Ga\bs G$ be the Cartan involutive homeomorphism
$\Ga g\mapsto \Ga{\;}^tg^{-1}$, so that for every bounded continuous
function $f:\RR/\ZZ \times\Ga\bs G\to\RR$ and for every $s\in\RR$, we
have
\[\int f \;dx\otimes d\,\Theta_*\,(\Phi^{-s})_*\,\mu_{\Ga_H\bs H}=
\int_{(x,y)\in(\RR/\ZZ)\times(\Ga_H\bs H)}f(x,\Theta(y\Phi^{-s}))
\,dx\,d\mu_{\Ga_H\bs H}(y)\;.\]

\bcoro [{{\bf Marklof \cite[Thm.~6]{Marklof10Inv}}}]\label{coro:marklof}
For every $t_0\in\RR$, for the narrow convergence of measures on
$\RR/\ZZ \times\Ga\bs G$, we have
\begin{equation}\label{eq:incoromarklof}
\lim_{t\to+\infty}\;\;\frac 1{\card\,\F_{t-t_0}}
\sum_{r\in\F_{t-t_0}}\Delta_r\otimes\Delta_{\Ga n_-(r)\Phi^t}= e^{t_0}
\int_{s=t_0}^{+\infty}
dx\otimes d\,\Theta_*\,(\Phi^{-s})_*\,\mu_{\Ga_H\bs H}\;e^{-s}\,ds\,.
\end{equation}
\ecoro

\dem We consider in this proof $X=\hdr$, where $\hnr$ is the upper
halfspace model of the real hyperbolic space of dimension $n$ (with
constant sectional curvature $-1$). We again denote by
$\iota:T^1\hnr\ra T^1\hnr$ the antipodal map $v\mapsto -v$. We
normalise, as we may, the Patterson density $(\mu_x)_{x\in X}$ of the
(nonuniform arithmetic) lattice $\Ga$ of the orientation preserving
isometry group $G$ of $X$ to consist of probability measures. The
critical exponent of $\Ga$ is
\begin{equation}\label{eq:deltaPSL2Z}
\delta_\Ga=1\;.
\end{equation}

We start the proof by recalling precisely a bijection between $G$ and
the unit tangent bundle of $\hdr$.  We denote by $\cdot$ the action of
$G$ by homographies on $\hdr\cup \partial_\infty \hdr$, as well at its
derived action on $T^1\hdr$. We fix $v^\bullet=(i,-i)\in T^1\hdr$,
which is the unit tangent vector at the base point $i$ of $\hdr$
pointing vertically down (its length is not adequate in the picture
below, but this makes the picture easier to understand).

\begin{center}
  \begin{picture}(0,0)%
\includegraphics{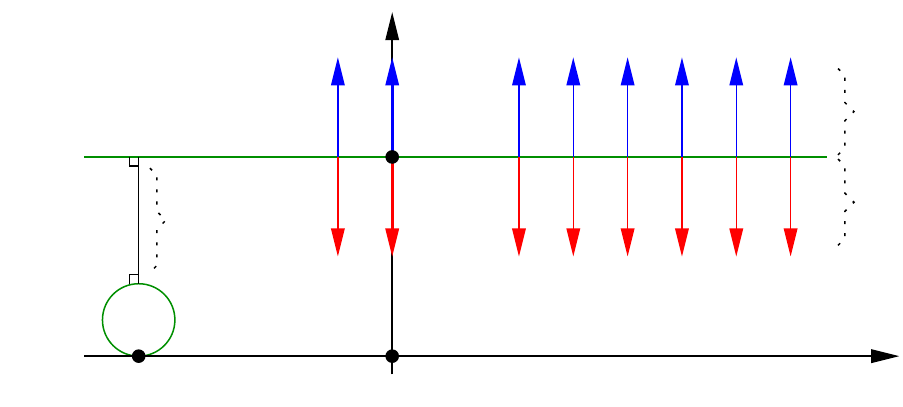}%
\end{picture}%
\setlength{\unitlength}{3812sp}%
\begingroup\makeatletter\ifx\SetFigFont\undefined%
\gdef\SetFigFont#1#2#3#4#5{%
  \reset@font\fontsize{#1}{#2pt}%
  \fontfamily{#3}\fontseries{#4}\fontshape{#5}%
  \selectfont}%
\fi\endgroup%
\begin{picture}(4485,2004)(121,-2740)
\put(991,-1861){\makebox(0,0)[lb]{\smash{{\SetFigFont{11}{13.2}{\rmdefault}{\mddefault}{\updefault}{\color[rgb]{0,0,0}$2\ln q$}%
}}}}
\put(2116,-871){\makebox(0,0)[lb]{\smash{{\SetFigFont{11}{13.2}{\rmdefault}{\mddefault}{\updefault}{\color[rgb]{0,0,0}$i\RR$}%
}}}}
\put(4591,-2536){\makebox(0,0)[lb]{\smash{{\SetFigFont{11}{13.2}{\rmdefault}{\mddefault}{\updefault}{\color[rgb]{0,0,0}$\RR$}%
}}}}
\put(1936,-1456){\makebox(0,0)[lb]{\smash{{\SetFigFont{11}{13.2}{\rmdefault}{\mddefault}{\updefault}{\color[rgb]{0,0,0}$i$}%
}}}}
\put(406,-1411){\makebox(0,0)[lb]{\smash{{\SetFigFont{11}{13.2}{\rmdefault}{\mddefault}{\updefault}{\color[rgb]{0,.56,0}$\H_\infty$}%
}}}}
\put(2116,-1861){\makebox(0,0)[lb]{\smash{{\SetFigFont{11}{13.2}{\rmdefault}{\mddefault}{\updefault}{\color[rgb]{1,0,0}$v^\bullet$}%
}}}}
\put(2116,-1231){\makebox(0,0)[lb]{\smash{{\SetFigFont{11}{13.2}{\rmdefault}{\mddefault}{\updefault}{\color[rgb]{0,0,1}$-v^\bullet$}%
}}}}
\put(4411,-1321){\makebox(0,0)[lb]{\smash{{\SetFigFont{11}{13.2}{\rmdefault}{\mddefault}{\updefault}{\color[rgb]{0,0,1}$W^+(-v^\bullet)$}%
}}}}
\put(4411,-1771){\makebox(0,0)[lb]{\smash{{\SetFigFont{11}{13.2}{\rmdefault}{\mddefault}{\updefault}{\color[rgb]{1,0,0}$W^-(v^\bullet)$}%
}}}}
\put(136,-2311){\makebox(0,0)[lb]{\smash{{\SetFigFont{11}{13.2}{\rmdefault}{\mddefault}{\updefault}{\color[rgb]{0,.56,0}$\ga\H_\infty$}%
}}}}
\put(541,-2671){\makebox(0,0)[lb]{\smash{{\SetFigFont{11}{13.2}{\rmdefault}{\mddefault}{\updefault}{\color[rgb]{0,0,0}$\frac{p}{q}=\ga\cdot\infty$}%
}}}}
\end{picture}%

\end{center}

We denote by $\wt \varphi:G\ra T^1\hdr$ the orbital map at
$v^\bullet$, defined by $g\mapsto g\cdot v^\bullet$, which is a
$G$-equivariant (for the left actions) homeomorphism, and by $\varphi:
\Ga\bs G\ra \Ga\bs T^1\hdr$ its quotient homeomorphism.  We define
$S=\begin{bmatrix} 0 & -1\\1&0\end{bmatrix}$, which is an order $2$
element of $\Ga$.  The involution $S$ satisfies the following
remarkable properties, in the connected centerfree semisimple real
Lie group $G$, that it anti-commutes with the standard Cartan sugbroup
$\Phi^\RR= \{\Phi^t:t\in\RR\}$ of $G$ and that the conjugation by $S$
is the standard Cartan involution $g\mapsto\;^tg^{-1}$ of $G$~:
\begin{equation}\label{eq:geodinvol}
\forall\;g\in G, \quad ^tg^{-1}=SgS^{-1}\qquad {\rm and}\qquad
\forall\;s\in \RR, \quad S\Phi^sS^{-1}=\Phi^{-s}\;.
\end{equation}
Hence, with $\Theta$ defined just before the statement of Corollary
\ref{coro:marklof},  for all $x\in \Ga\bs G$ and
$s\in\RR$, we have
\[
\Theta(x\Phi^s)=\Theta(x)\Phi^{-s}\;.
\]
The element $S$ represents a generator of the order $2$ standard Weyl
group $N_G(\Phi^\RR)/Z_G(\Phi^\RR)$. The following properties say that
the action of the geodesic flow $\flow{t}$ on $T^1\hdr$ corresponds to
the multiplication on the right by $\Phi^t$ in $G$, and that the
antipodal map on $T^1\hdr$ corresponds to the multiplication on the
right by $S$ in $G$~:
\begin{equation}\label{eq:conjugflogeod}
  \forall\;t\in \RR,\;\forall\; g\in G,\quad
  \flow{t}\wt\varphi(g)=\wt\varphi(g\Phi^t)
  \quad{\rm and}\quad \iota\;\wt\varphi(g)=\wt\varphi(gS)\;.
\end{equation}
By the above two centered formulas and since $S\in\Ga$, the
homeomorphism $\varphi$ relates the antipodal map $\iota$ on $\Ga\bs
T^1\hdr$ to the Cartan involution $\Theta$ on $\Ga\bs G$ by
\[
\iota\circ\varphi =\varphi\circ\Theta\;.
\]

Let $\H_\infty=\{z\in \hdr:\Im \;z\geq 1\}$, which is a (closed)
horoball centered at $\infty$ in $\hdr$. The subgroup $\Ga_H$ is equal
to the stabiliser $\Ga\!_{\H_\infty}$ of $\H_\infty$ in $\Ga$. We define
\begin{equation}\label{eq:defDpm}
\D^-=\D^+=(\ga\H_\infty)_{\ga\in\Ga}\;,
\end{equation}
which are locally finite $\Ga$-equivariant families of horoballs. The
map from $\Ga-\Ga\!_{\H_\infty}$ to $\RR$ defined by $\ga=
\begin{bmatrix} p & r\\q&s\end{bmatrix}\mapsto \ga\cdot\infty
=\frac{p}{q}$ (where we assume, as we may, that $q>0$) induces a
bijection from $\Ga\!_{\H_\infty}\bs(\Ga-\Ga\!_{\H_\infty})/
\Ga\!_{\H_\infty}$ to the additive group $\QQ/\ZZ$ such that
$d(\H_\infty,\ga\H_\infty) =2\ln q$ (see the above picture).  In
particular, for all $t,t_0\in\RR$, we have
\[
d(\H_\infty,\ga\H_\infty)\leq t-t_0\quad{\rm if~and~only~if}\quad
q\leq e^{\frac{t-t_0}{2}}\;.
\]
Identifying geodesic rays in $\G_{+,0}X$ and geodesic lines in $\G X$
with their unit tangent vector at time $0$, we have
\[
\normalout \H_\infty=W^-(v^\bullet)=\wt \varphi(H)\;,
\]
so that, by the left equivariance of $\wt\varphi$, the orbits of the
right action of $H$ on $G$ correspond to the strong unstable leaves
for the geodesic flow on $T^1\hdr$.  Similarly, using Equation
\eqref{eq:conjugflogeod}, we have
\[
\normalin\H_\infty=W^+(-v^\bullet)=\iota W^-(v^\bullet)=
\wt \varphi(H\,S)\quad{\rm and}\quad
W^{0+}(-v^\bullet)=\wt \varphi(H\Phi^\RR\,S)\;.
\]
More precisely, using the right part of Equation \eqref{eq:geodinvol}
and Equation \eqref{eq:conjugflogeod}, we have
\begin{equation}\label{eq:weakstable}
  \forall\;s,r\in\RR,\quad \wt \varphi(n_-(r)\,\Phi^{-s}\,S)=
  \wt \varphi(n_-(r)\,S\,\Phi^{s})=
  \flow{s} \;\iota\,\wt \varphi(n_-(r))\;.
\end{equation}

The endpoint map $\wt\psi_\pm:\normalpm\H_\infty\ra\RR$ defined by
$\rho\mapsto \rho_\pm$ is a $\Ga_H$-equivariant homeomorphism, such that
$\wt \varphi^{-1}(\wt\psi_+^{-1}(r))=n_-(r)$ for all $r\in\RR$. We
denote by $\psi_\pm:\Ga_H\bs\normalpm\H_\infty\ra\RR/\ZZ$ the
quotient homeomorphism, and we identify $\Ga_H\bs\normalout
\H_\infty$ with its image in $\Ga\bs T^1\hdr$. For every $\ga
=\begin{bmatrix} p & r\\q&s \end{bmatrix} \in \Ga-\Ga_H$, with
$\rho_\ga\in\normalout\H_\infty$ the geodesic ray entering
perpendicularly in $\ga \H_\infty$, we have
\begin{equation}\label{eq:rhogamprop}
\wt\varphi^{-1}(\rho_\ga)=n_-(\ga\cdot\infty)
\quad{\rm and}\quad
(\psi_+)_*(\Delta_{\Ga\rho_\ga})=\Delta_{\ga\cdot\infty\!\!\!\!\mod 1}=
\Delta_{\frac{p}{q}\!\!\!\!\mod 1}\;.
\end{equation}
Furthermore, by \cite[Thm.~9.11]{ParPau16LMS} or \cite[Prop.~20
  (2)]{ParPau17ETDS} with $n=2$, the skinning measure
$\wt\sigma^\pm_{\H_\infty}$ is equal to twice the Riemannian volume of
$\normalmp \H_\infty$, so that
\begin{equation}\label{eq:endpointskinn}
  (\psi_+)_*(\sigma^+_{\D^-})=2\,dx\quad{\rm and}\quad \
  (\varphi^{-1})_*(\sigma^+_{\D^-})=2\,\mu_{\Ga_H\bs H}\;.
\end{equation}

By for instance \cite[Thm.~9.10]{ParPau16LMS} or \cite[Prop.~20
  (1)]{ParPau17ETDS} with $n=2$, we have
\[
\|m_{\rm BM}\|=4\pi\vol(\Ga\bs\hdr)=\frac{4\pi^2}3\;.
\]
Mertens's formula \cite[Thm.~330]{HarWri08} (see also \cite[\S 3]
{ParPau14AFST} for a geometric proof) implies that, as $t\ra+\infty$,
\[
\card\,\F_{t-t_0}\sim\frac{3}{\pi^2}e^{2\frac{t-t_0}{2}}=
\frac{3}{\pi^2}e^{t-t_0}\;.
\]
Since no element of $\Ga$ pointwise fixes a nontrivial geodesic
segment of $\hdr$, for every $\ga\in\Ga$ such that
$d(\H_\infty,\ga\H_\infty)>0$, we have
\[
m_\ga=1\;.
\]

For every $t_0\in\RR$, let us consider the truncation $\Phi^{\geq t_0}
=\{\Phi^t:t\geq t_0\}$ of the Cartan subgroup $\Phi^\RR$. For all
$t\in\RR$ and $\ga\in\Ga-\Ga_H$, by the two left parts of Equations
\eqref{eq:conjugflogeod} and \eqref{eq:rhogamprop}, we have
\[
(\varphi^{-1})_*(\Delta_{\Ga\flow{t}\rho_\ga})=
\Delta_{\Ga n_{-}(\ga\cdot\infty)\Phi^{t}}\;.
\]
By Equation \eqref{eq:weakstable}, the homeomorphism $\varphi^{-1}$
maps the truncated stable leaf
\[
\Ga W^{0+}_{t_0}(-v^\bullet)= \bigcup_{s\geq t_0}
\Ga\flow{s}\partial^1_-\H_\infty = \bigcup_{s\geq t_0}
\Ga\flow{s}W^{+}(-v^\bullet) =\bigcup_{s\geq t_0} \Ga \flow{s} \iota
W^{-}(v^\bullet)
\]
to the truncated orbit $\Ga H(\Phi^{\geq t_0})^{-1}S$ in $\Ga\bs G$ of
the lower triangular subgroup of $G$. Besides, by the left part of
Equation \eqref{eq:geodinvol} and since $S\in \Ga$ for the first
equality, by Equation \eqref{eq:weakstable} for the third equality, by
Equations \eqref{eq:strongstabmeasj} with \eqref{eq:deltaPSL2Z} for
the fourth equality, and since
$\iota_*\sigma^-_{\D^+}=\sigma^+_{\D^-}$ and by the right part of
Equation \eqref{eq:endpointskinn} for the last equality, for all
$s,r\in\RR$ with $s\geq t_0$, we have
\begin{align*}
  &d\big((\varphi^{-1})_*(\mu^{0+}_{\D^+,t_0})\big)
  (\Theta(\Ga n_{-}(r)\Phi^{-s})) =
  d\big((\varphi^{-1})_*(\mu^{0+}_{\D^+,t_0})\big)
  (\Ga n_{-}(r)\Phi^{-s}S)\\=\;&
  d\mu^{0+}_{\D^+,t_0}\big(\Ga\wt\varphi(n_{-}(r)\Phi^{-s}S)\big) =
  d\mu^{0+}_{\D^+}\big(\Ga\flow{s}\iota\, \wt\varphi(n_{-}(r))\big) =
  e^{-s}d\sigma^-_{\D^+}\big(\Ga\iota\, \wt\varphi(n_{-}(r))\big)\,ds
  \\ =\;&
  e^{-s}(\varphi^{-1})_*\,\iota_*\,d\sigma^-_{\D^+}(\Ga n_{-}(r))\,ds
  =2\;d\mu_{\Ga_H\bs H}(\Ga n_{-}(r))\,e^{-s}\,ds\;.
\end{align*}
Therefore, by the left part of Equation \eqref{eq:endpointskinn}, for
all $x\in\RR/\ZZ$, $y\in\Ga_H\bs H$ and $s\geq t_0$, we have
\begin{equation}\label{eq:leftmarklof}
  d\big( (\psi_+\times \varphi^{-1})_*\big(\sigma^+_{\D^-}\otimes
  \mu^{0+}_{\D^+,t_0}\big)\big)(x,\Theta(y\,\Phi^{-s}))=
  4\;dx\;d\mu_{\Ga_H\bs H}(y)\,e^{-s}\,ds\;.
\end{equation}
By the linearity of the pushforward of measures and by Equation
\eqref{eq:rhogamprop}, for every $t\in\RR$, we hence have
\begin{align}
 (\psi_+\times \varphi^{-1})_*&\Big(\|m_{\rm BM}\|\;e^{-\delta_\Ga\, t}
\sum_{\substack{\ga\,\in\;\Ga\!_{\H_\infty}\!\backslash \Ga/\Ga\!_{\H_\infty})\\ 
0<d(\H_\infty,\ga \H_\infty)\leq\, t-t_0}} \; 
m_\ga \;\Delta_{\Ga\rho_\ga}\otimes\Delta_{\flow t\Ga\rho_\ga}\Big)
\nonumber\\&=\frac{4\pi^2}3\;e^{-t}\sum_{r\in\F_{t-t_0}}
\Delta_{r}\otimes\Delta_{\Ga n_-(r)\Phi^t}\nonumber\\
&\sim
4\,e^{-t_0}\;\frac 1{\card\,\F_{t-t_0}}\sum_{r\in\F_{t-t_0}}
\Delta_{r}\otimes\Delta_{\Ga n_-(r)\Phi^t}\;.
\label{eq:rightmarklof}\end{align}
Since the product map $\psi_+\times \varphi^{-1}$ is a homeomorphism
from $(\Ga W^-(v^\bullet)) \times (\Ga W^{0+}_{t_0}(-v^\bullet))$ to
$(\RR/\ZZ) \times (\Ga H(\Phi^{\geq t_0})^{-1})$, its pushforward map
on measures is continuous for the narrow convergence. Hence Corollary
\ref{coro:marklof} follows from Equations \eqref{eq:leftmarklof} and
\eqref{eq:rightmarklof} by Theorem \ref{theo:quotient} (1) applied to
the families $\D^\pm$ defined in Equation \eqref{eq:defDpm}.  \cqfd

\medskip\noindent {\bf Remarks. } (1) Using the final Remark of
Section \ref{sec:shrinking} and an approximation by linear
combinations of functions with separate variables, one could prove
that there exist $\tau'>0$ and $\ell\in\NN$ such that for every
$\psi'\in \C^\ell_c(\RR/\ZZ\times\Ga\bs G)$, we have an error term of
the form $\bigO_{t_0}(e^{-\kappa' t}\|\psi'\|_\ell)$ when replacing
(before taking the limit on the left hand side) in the two sides of
Equation \eqref{eq:incoromarklof} the function $f$ by $\psi'$. See
also \cite{Marklof13} when $n=2$ and \cite{Li15} when $n\geq 3$ for an
effective version of Marklof's result.

\medskip
(2) A version of Corollary \ref{coro:marklof} with congruences is
possible.  Let $N\in \NN-\{0\}$, and let $\Ga_0[N]$ be the Hecke
congruence subgroup of level $N$ of $\Ga$, preimage of the upper
triangular subgroup by the morphism of reduction modulo $N$ of the
coefficients. Up to replacing $\F_t$ by $\{\frac{p}{q}\!\mod 1:
\frac{p}{q}\in \F_t,\; q\equiv 0\!\mod N\}$, to replacing $\Ga$ by
$\Ga_0[N]$ and to replacing $\Theta_*$ by an averaging operator over
cosets of $\Ga_0[N]$ in $\Ga$ (coming from the fact that the lattice
$\Ga_0[N]$ is no longer invariant under the Cartan involution
$g\mapsto \;^tg^{-1}$), one could obtain as in \cite[Thm.~2
  (B)]{Marklof10AMS} a joint partial equidistribution of Farey
fractions with a congruence assumption on their denominator and with
an error term. See also \cite{Heersink21}.

\subsection{Equidistribution of complex Farey fractions
  at a given density}
\label{subsec:complexFarey}

Let $K$ be an imaginary quadratic number field, with discriminant
$D_K$, ring of integers $\OOO_K$, finite group of unit integers
$\OOO_K^\times$ (which is equal to $\{\pm 1\}$ unless $D_K=-4,-3$),
and Dedekind's zeta function $\zeta_K$.

Let $G=\PSLC$ and let $\Ga$ be the {\it Bianchi group}
$\PSL_2(\OOO_K)$.  For all $r\in\CC$ and $t\in \RR$, we consider the
elements of $G$ defined by
\[
n_-(r)=\begin{bmatrix} 1 & r\\0&1\end{bmatrix}\quad 
{\rm and}
\quad 
\Phi^t=\begin{bmatrix} \;e^{-t/2} & 0\\0&e^{t/2}\end{bmatrix}\,.
\]
Let $H=\{n_-(r):r\in\CC\}$. We denote by $$M=\Big\{\begin{bmatrix}
\;e^{-i\,\theta/2} & 0\\0&e^{i\,\theta/2}\end{bmatrix}: \theta\in
\RR\Big\}$$ the compact factor of the centraliser of the standard
Cartan subgroup $\Phi^\RR=\{\Phi^t:t\in\RR\}$ of $G$, which normalises
$H$. Note that both $\Ga$ and $M$ are invariant under the standard
Cartan involution $g\mapsto \;^tg^{-1}$.  Let 
\[
\Ga_H=N_G(H)\cap\Ga=
(HM)\cap \Ga= \Big\{\begin{bmatrix} a & b\\0&a^{-1}\end{bmatrix}:
a\in\OOO_K^\times,\;b\in\OOO_K\Big\}\,,
\]
which is a semi-direct product $(M\cap\Ga)\semidirect (H\cap \Ga)$.
The discrete group $\Ga_H$ acts properly discontinuously on the left
on $H$ so that $H\cap \Ga$ acts firstly by translations and $M\cap
\Ga$ secondly by conjugation: for all $a\in\OOO_K^\times$,
$b\in\OOO_K$ and $r\in\CC$, we have
\[
\begin{bmatrix} a & b\\0&a^{-1}\end{bmatrix}\cdot \begin{bmatrix} 1 &
  r\\0&1\end{bmatrix} = \begin{bmatrix} 1 & a^2r+ab\\0&1\end{bmatrix}\,.
\]
We see, as we may, $\Ga_H\bs H$ contained in $\Ga\bs G/M$ (as the
image of $H$ in the set of double cosets). We endow $\Ga_H\bs H$ with
the induced measure $\mu_{\Ga_H\bs H}$ of a Haar measure on $H$ by the
branched cover $H\ra \Ga_H\bs H$, normalised to be a probability
measure, that we also see as a probability measure on $\Ga\bs G/M$
(with support $\Ga_H\bs H=(\Ga\cap H)\bs H/(H\cap M)$). We denote by
$\OOO'_K$ the semidirect product $\OOO_K^\times \semidirect \OOO_K$,
which acts on the left, with kernel of order $2$, on $\CC$ by
$((a,b),r)\mapsto a^2r+ab$.  Note that for every $t\in\RR$, since
$\Phi^t$ centralises $M$, the double class $\Ga n_-(r)\Phi^tM$ is well
defined for every equivalence class $r\in\OOO'_K\bs\CC$. We endow the
quotient space $\OOO'_K\bs \CC$ with the induced measure $dx$ of the
Lebesgue measure on $\CC$ by the branched cover $\CC\mapsto \OOO'_K\bs
\CC$, normalised to be a probability measure.

For every $t\in\RR$, we consider the subset $\F_{t}$ of
$\OOO'_K\bs\CC$ consisting of the {\em complex Farey fractions of
  height at most} $e^{t/2}$, defined by
\[
\F_{t}=\OOO'_K\bs\Big\{\frac pq:\;p,q\in\OOO_K,
\quad p\OOO_K+q\OOO_K=\OOO_K,\quad 0<|q|\leq e^{t/2}\Big\}\,.
\]
Note that the above set of fractions $\frac pq$ is indeed invariant
under $\OOO'_K$.

Let $\Theta:\Ga\bs G/M\ra \Ga\bs G/M$ be the Cartan
involutive homeomorphism defined by $\Ga gM\mapsto \Ga{\;}^tg^{-1}M$.

\bcoro\label{coro:complexFarey} For every $t_0\in\RR$, for the narrow
convergence of probability measures on $(\OOO'_K\bs \CC) \times(\Ga\bs
G/M)$, we have
\begin{align*}
\lim_{t\to+\infty}\;\;\frac 1{\card\,\F_{t-t_0}}&
\sum_{r\in\F_{t-t_0}}\Delta_r\otimes\Delta_{\Ga n_-(r)\Phi^tM}\\ &= 2\;e^{2t_0}
\int_{s=t_0}^{+\infty}(dx)\otimes\big(
\Theta_*\,(\Phi^{-s})_*\,\mu_{\Ga_H\bs H}\big)\,e^{-2s}\,ds\,.
\end{align*}
\ecoro

This statement implies Corollary \ref{coro:complexFareyintro} in the
introduction when $D_K\neq -4,-3$, since then $\OOO'_K\bs\CC=
\OOO_K\bs\CC=\CC/\OOO_K$ and $\Ga_H=H\cap \Ga$. As a remark similar to
the remark at the end of Section \ref{subsec:retrouvonsMarklof}, one
could obtain an error term under an additional regularity assumption, and
a joint partial equidistribution result of complex Farey fractions
with their denominator congruent to $0$ modulo any fixed element $N$
in $\OOO_K-\{0\}$.

\medskip
\dem We mostly indicate the differences with the proof of Corollary
\ref{coro:marklof}. We now consider $X=\htr$ with coordinates
$(z,u)\in\CC\times\; ]0,+\infty[$. The critical exponent of the
  (nonuniform arithmetic) lattice $\Ga$ of the orientation preserving
  isometry group $G$ of $X$ is now
\[
\delta_\Ga=2\;.
\]
We denote by $\cdot$ the action of $G$ by homographies on $
\partial_\infty \htr=\CC\cup\{\infty\}$, by isometries on $\htr$
through the Poincaré extension, and by the derived action on
$T^1\htr$. We now fix $v^\bullet=((0,1),(0,-1))\in T^1\htr$. Its
stabiliser in $G$ is equal to $M$ and is hence centralised by
$\Phi^\RR$.  The orbital map $\wt\varphi:g\mapsto g\cdot v^\bullet$
now defines a homeomorphism $\varphi: \Ga\bs G/M\ra \Ga\bs T^1\htr$.
The order $2$ element $S=\begin{bmatrix} 0 & -1\\1&0\end{bmatrix}$
still belongs to $\Ga$. It normalises $M$ and $\Phi^\RR$, and the
formulae \eqref{eq:geodinvol} and \eqref{eq:conjugflogeod} are still
satisfied.

Let now $\H_\infty=\{(z,u)\in \htr:u\geq 1\}$. With $\Phi^{\geq t_0}
=\{\Phi^t:t\geq t_0\}$, we again have
\begin{equation}\label{eq:relatnormalinoutwtphi}
\normalout\H_\infty=W^-(v^\bullet)=\wt \varphi(H)\quad{\rm and}\quad
W^{0+}_{t_0}(-v^\bullet)
= \bigcup_{s\geq t_0}
\flow{s}\normalin\H_\infty=\wt \varphi(H(\Phi^{\geq t_0})^{-1}S)\;.
\end{equation}
The subgroup $\Ga_H$ is again equal to the stabiliser
$\Ga\!_{\H_\infty}$ of the horoball $\H_\infty$ in $\Ga$. We again
consider the locally finite $\Ga$-equivariant families of horoballs
\[
\D^+=\D^-=(\ga\H_\infty)_{\ga\in\Ga}\;.
\]
The map $\ga= \begin{bmatrix} p & r\\q&s\end{bmatrix}\mapsto
\ga\cdot\infty =\frac{p}{q}$ now induces, for every $t\in\RR$, a
bijection from $\big\{[\ga]\in\Ga\!_{\H_\infty}\bs\Ga/\Ga\!_{\H_\infty}:
0<d(\H_\infty, \ga\H_\infty)\leq t\big\}$ to $\F_t$.  With $\rho_\ga$
the element of $\normalout\H_\infty$ whose point at infinity is
$\ga\cdot\infty$, the endpoint map $\wt\psi_\pm: \normalpm
\H_\infty\ra\CC$ now induces an homeomorphism $\psi_\pm: \Ga_H\bs
\normalpm\H_\infty \ra\OOO'_K\bs \CC$, such that
\[
(\psi_+)_*(\Delta_{\Ga\rho_\ga})=\Delta_{\OOO'_K\ga\cdot\infty}\;.
\]

Let us compute the total mass of the induced Lebesgue measure
$d\operatorname{Leb}_{\OOO'_K \bs \CC}$ on $\OOO'_K\bs \CC$, yielding
$dx$ after renormalisation to a probability measure. Since the
branched cover $\OOO_K\bs\CC\ra\OOO'_K\bs\CC$ is
$\frac{|\OOO_K^\times|}{2}$-sheeted outside the singular part and
since $\OOO_K$ is generated as a $\ZZ$-lattice of $\CC$ by $1$ and
$(D_K+i\sqrt{|D_K|})/2$, we have
\[
\|d\operatorname{Leb}_{\OOO'_K\bs \CC}\|=\frac{2}{|\OOO_K^\times|}\;
\|d\operatorname{Leb}_{\OOO_K\bs \CC}\|=
\frac{\sqrt{|D_K|}}{|\OOO_K^\times|}\;.
\]
Again by \cite[Thm.~9.11]{ParPau16LMS} or \cite[Prop.~20
  (2)]{ParPau17ETDS} with now $n=3$, we have
\[
(\psi_+)_*(\sigma^+_{\D^-})=4\,d\operatorname{Leb}_{\OOO'_K\bs \CC}=
\frac{4\,\sqrt{|D_K|}}{|\OOO_K^\times|}\,dx
\quad {\rm and}\quad (\varphi^{-1})_*(\sigma^+_{\D^-})
=\frac{4\,\sqrt{|D_K|}}{|\OOO_K^\times|}\,d\mu_{\Ga_H\bs H}\;.
\]
Again by \cite[Thm.~9.10]{ParPau16LMS} or \cite[Prop.~20
  (1)]{ParPau17ETDS} with now $n=3$ and with Humbert's volume formula (see
for instance \cite[\S 8.8 and \S 9.6]{ElsGruMen98}), we have
\[
\|m_{\rm BM}\|=4\Vol(\SSS^2)\Vol(\Ga\bs\htr)
=\frac{4}{\pi}\;|D_K|^{3/2}\;\zeta_K(2)\;.
\]
Mertens's formula for the quadratic imaginary fields (see also
\cite[Theo.~3.1] {ParPau14AFST}) gives, using the action of
$k\in\OOO_K$ on $(p,q)\in\OOO_K\times\OOO_K$ by horizontal shears
$k\cdot (p,q)=(p+kq)$, as $t\ra+\infty$,
\begin{align*}
  \card\,\F_{t-t_0}&\sim\frac{2}{|\OOO_K^\times|}
  \card\big(\OOO_K\bs\Big\{\frac pq:\;p,q\in\OOO_K,
  \; p\OOO_K+q\OOO_K=\OOO_K,\; 0<|q|\leq e^{(t-t_0)/2}\Big\}\big)
  \\&=\frac{2}{|\OOO_K^\times|^2} \card\big(\OOO_K\bs\{(p,q)\in
\OOO_K\times\OOO_K:p\OOO_K+q\OOO_K=\OOO_K, \;0<|q|^2\leq e^{t-t_0}\}\big)
\\&\sim\frac{2\,\pi}{|\OOO_K^\times|^2\;\zeta_K(2)\sqrt{|D_K|}}\;
  e^{2t-2t_0}\;.
\end{align*}
Since $\OOO_K$ has finite index in $\OOO'_K$, there are only finitely
many elliptic elements in $\Ga$ up to conjugation by $\Ga\cap H$ whose
fixed point set contains $\infty$ as a point at infinity. There are
only finitely many $\Ga\!_{\H_\infty}$-orbits of images of $\H_\infty$
by $\Ga$ meeting $\H_\infty$. Hence there exists a finite subset $F$
of the set of double cosets $\Ga\!_{\H_\infty}\bs\Ga/ \Ga\!_{\H_\infty}$
such that for every element $[\ga]\in\Ga\!_{\H_\infty}\bs\Ga/
\Ga\!_{\H_\infty}-F$, we have
\[
d(\H_\infty,\ga\H_\infty)>0 \quad {\rm and}\quad m_\ga=1\;.
\]
We have similarly, for all $\ga\in \Ga-\Ga\!_{\H_\infty}$ and $t\in\RR$,
\[
(\varphi^{-1})_*(\Delta_{\Ga\flow{t}\rho_\ga})=
\Delta_{\Ga n_{-}(\ga\cdot\infty)\Phi^{t}M}
\]
and, for all $y\in \Ga_H\bs H$ and $s\in\RR$ with $s\geq t_0$,
\begin{align*}
d\big((\varphi^{-1})_*(\mu^{0+}_{\D^+,t_0})\big)
(\Theta(y\,\Phi^{-s})) &=\|\sigma^-_{\D^+}\|\;
d\mu_{\Ga_H\bs H}(y)\,e^{-2s}\,ds\\ &
=\frac{4\;\sqrt{|D_K|}}{|\OOO_K^\times|}\;
d\mu_{\Ga_H\bs H}(y)\,e^{-2s}\,ds\;.
\end{align*}
The end of the proof of Corollary \ref{coro:complexFarey} proceeds
now as the one of Corollary \ref{coro:marklof}.
\cqfd

\subsection{Equidistribution of Heisenberg Farey fractions
  at a given density}
\label{subsec:heisenbergFarey}

Let $K,D_K,\OOO_K,\OOO_K^\times,\zeta_K$ be as in the beginning of
Section \ref{subsec:complexFarey}. Let $\tr$ and $\n$ be the (absolute)
trace and norm of $K$. We denote by $\langle a, \alpha,c\rangle$ the
ideal of $\OOO_K$ generated by $a,\alpha,c\in\OOO_K$.

Let $q$ be the nondegenerate Hermitian form $-z_0\overline{z_2}
-z_2\overline{z_0}+|z_1|^2$ of signature $(1,2)$ on $\CC^3$ with
coordinates $(z_0,z_1,z_2)$. Let $G=\PSU_q=\SU_q/(\UU_3\id)$ be the
projective special unitary group of $q$, where $\SU_q=\{g\in
\GL_3(\CC): q\circ g=q,\;\det g=1\}$ and $\UU_3$ is the group of cube
roots of unity. Let $\Ga$ be the image of $\SU_q\cap\SL_3(\OOO_K)$ in
$G$, which is a (nonuniform) arithmetic lattice in $G$, called the
(projective special) {\em Picard modular group} of $K$.

Denoting by $\begin{bmatrix} a & \overline{\ga} & b \\
  \alpha & A & \beta\\c & \overline{\delta} & d \end{bmatrix}$
the image in $G$ of $\begin{pmatrix} a & \overline{\ga} & b \\
  \alpha & A & \beta\\c & \overline{\delta} & d \end{pmatrix}\in
\SU_q$, let
\[
H=\Big\{n_-(w_0,w)=\begin{bmatrix} 1 & \overline{w} & w_0\\0&1&w
\\0&0&1\end{bmatrix}:w_0,w\in\CC, \;2\;\Re \; w_0=|w|^2\Big\}\,,
\]
\[
\Phi^\RR=\Big\{
\Phi^t=\begin{bmatrix} \;e^{-t} & 0& 0\\0&1&0\\0&0&e^{t}\end{bmatrix}
:t\in\RR\Big\}\quad {\rm and}\quad  M=\Big\{
\begin{bmatrix} \;\zeta & 0& 0\\0&\overline\zeta^{\,2}&0\\0&0
  &\zeta\end{bmatrix}:\zeta\in\CC,\;|\zeta|=1\Big\}\,.
\]
Note that $H$, $\Phi^\RR$ and $M$ are Lie subgroups of $G$, that $M$
is the compact factor of the centraliser in $G$ of the standard Cartan
subgroup $\Phi^\RR$ of $G$, and that the subgroup $M\Phi^\RR$
normalises the {\em Heisenberg group} $H$. The groups $\Ga$ and $M$
are invariant under the standard Cartan involution $$g\mapsto
\;^*g^{-1}\,,$$ where $^*g$ is the image in $G$ of the
transpose-conjugate matrix of any matrix in $\SU_q$ representing $g$.

Let $\Ga_H=N_G(H)\cap\Ga= (MH)\cap \Ga= \Big\{\begin{bmatrix} \;u &
u\, \overline{v}& u\,v_0\\0&\overline u^{\,2}&\overline
u^{\,2}\,v\\0&0 &u\end{bmatrix}:\begin{array}{c} u\in\OOO_K^\times,
  \;v,v_0\in\OOO_K\\ \tr(v_0)=\n(v)\end{array}\Big\}$, which acts
properly discontinuously on the left on $H$ by
\begin{equation}\label{eq:actionGasubHsurH}
  \begin{bmatrix} \;u & u\,\overline{v}& u\,v_0\\
    0&\overline u^{\,2}&\overline u^{\,2}\,v\\0&0&u
  \end{bmatrix}\cdot\begin{bmatrix} \;1 &
\overline{w}& w_0\\0&1&w\\0&0&1\end{bmatrix}
=\begin{bmatrix} \;1 &
u^3\,(\overline{w}+\overline{v})& w_0+v_0+w\,\overline{v}\\
0&1&\overline{u}^3\,(w+v)\\0&0&1\end{bmatrix}\,,
\end{equation}
where $H\cap \Ga$ acts firstly by left translations and $M\cap \Ga$
secondly by conjugations on $H$. The inclusion map $H\ra G$ induces an
identification between the quotient $\Ga_H\bs H$ and the image of $H$
in $\Ga\bs G/M$.  We endow $\Ga_H\bs H$ with the induced measure
$\mu_{\Ga_H\bs H}$ of a Haar measure on $H$, by the branched cover
$H\ra \Ga_H\bs H$, normalised to be a probability measure, that we
also see as a probability measure on $\Ga\bs G/M$ (with support
$\Ga_H\bs H$).

For every $t\in\RR$, we consider the subset $\F_{t}$ of $\Ga_H\bs H$
consisting of the {\em Heisenberg Farey fractions of height at most}
$e^{t}$, defined by
\[
\F_{t}=\Ga_H\bs\Big\{n_-\big(\frac{a}{c},\frac{\alpha}{c}\big):
\begin{array}{c}a,\alpha,c\in\OOO_K,\;\langle a,\alpha,c\rangle
  =\OOO_K,\\ \tr(a\,\overline{c})=\n(\alpha),\end{array}\quad
0<\n(c)\leq e^{2\,t}\Big\}\,.
\]
Note that the above set of elements $n_-\big(\frac{a}{c},
\frac{\alpha}{c}\big)$ is indeed invariant
under $\Ga_H$, by Equation \eqref{eq:actionGasubHsurH}.
Let $\Theta:\Ga\bs G/M\ra \Ga\bs G/M$ be the Cartan
involutive homeomorphism defined by $\Ga gM\mapsto \Ga{\;}^*g^{-1}M$.

\bcoro\label{coro:HeisenbergFarey} For every $t_0\in\RR$, for the
narrow convergence of probability measures on $(\Ga_H\bs H)
\times(\Ga\bs G/M)$, we have
\begin{align*}
\lim_{t\to+\infty}\;\;\frac 1{\card\,\F_{t-t_0}}&
\sum_{r\in\F_{t-t_0}}\Delta_r\otimes\Delta_{\Ga\, r\,\Phi^tM}\\&
=4\;e^{4\,t_0} \int_{s=t_0}^{+\infty}(\mu_{\Ga_H\bs H})\otimes
\big(\Theta_*\,(\Phi^{-s})_*\,\mu_{\Ga_H\bs H}\big)\,e^{-4s}\,ds\,.
\end{align*}
\ecoro

As a remark similar to the remark at the end of Section
\ref{subsec:retrouvonsMarklof}, one could obtain an error term under
an additional regularity assumption, and a joint partial
equidistribution result of Heisenberg Farey points
$n_-\big(\frac{a}{c}, \frac{\alpha}{c})$ modulo $\Ga_H$ with their
denominators $c$ congruent to $0$ modulo any fixed element $N$ in
$\OOO_K-\{0\}$.

\medskip
\dem We mostly indicate the differences with the proof of Corollary
\ref{coro:marklof}. We refer to \cite{Goldman99} as well as \cite[\S
  6.1]{ParPau10GT}, \cite[\S 3]{ParPau17MA} for background on complex
hyperbolic geometry. We follow the conventions of this last reference
concerning the normalisation of the sectional curvature and the choice
of the Hermitian form with signature $(1,2)$.

We now consider $X=\hdc$ the Siegel domain model of the complex
hyperbolic plane, that is, the complex manifold
\[
\big\{(w_0,w)\in\CC^2: 2\,\Re\; w_0 -|w|^2>0\big\}\,,
\]
endowed with the Riemannian metric
\begin{equation}\label{eq:metriemSiegeldom}
ds^2_{\,\hdc}=\frac{1}{(2\,\Re\; w_0 -|w|^2)^2}
\big((dw_0-dw\, \overline{w})(\overline{dw_0}-w\,\overline{dw})+
(2\,\Re\; w_0 -|w|^2)\;dw\, \overline{dw}\,\big)\,.
\end{equation}
This metric is normalised so that its sectional curvatures are in
$[-4,-1]$. The boundary at infinity of $\hdc$ is
\[
\partial_\infty\hdc=\big\{(w_0,w)\in
\CC^2 \;:\; 2\,\Re\; w_0 -|w|^2=0\big\}\cup\{\infty\}\,.
\]
Using homogeneous coordinates, we identify $\hdc\cup\partial_\infty
\hdc$ with its image in $\PP^2(\CC)$ by the map $(w_0,w)\mapsto
    [w_0:w:1]$ and $\infty \mapsto [1:0:0]$. We denote by $\cdot$ the
projective action of $G$ on $\hdc\cup\partial_\infty \hdc$, as well
as its derived action on $T^1\hdc$. The holomorphic isometry group of
$\hdc$ is $G$ (acting projectively on $\PP^2(\CC)$).

The critical exponent of the (nonuniform arithmetic) lattice $\Ga$ of
$G$ is now (see for instance \cite[\S 6]{CorIoz99})
\[
\delta_\Ga=4\;.
\]

We now fix $v^\bullet=((1,0),(-2,0))\in T^1\hdc$, which is indeed a
unit tangent vector with footpoint $x^\bullet =(1,0)$ by Equation
\eqref{eq:metriemSiegeldom}. The stabiliser of $v^\bullet$ in $G$ is
equal to $M$ and is hence centralised by $\Phi^\RR$. The orbital map
$\wt\varphi: g\mapsto g\cdot v^\bullet$ now defines a homeomorphism
$\varphi: \Ga\bs G/M \ra \Ga\bs T^1\hdc$.

For every $t\in\RR$, the element $\Phi^t$ acts on $\hdc$ by the map
$(w_0,w)\mapsto (e^{-2t}w_0,e^{-t}w)$. The geodesic line $\ell$ in
$\hdc$ such that $\ell(0)=x^\bullet$ and $\ell'(0)= v^\bullet$ is
$t\mapsto (e^{-2\,t},0)$. Hence $\flow{t}
v^\bullet=\ell'(t)=(-2\,e^{-2\,t},0)= d_x\Phi^t(v^\bullet) =\Phi^t
\cdot v^\bullet$, therefore by equivariance, we have
\[
\forall\;t\in\RR,\;\forall\;g\in G,\quad \flow{t}\wt\varphi(g)=
\wt\varphi(g\Phi^t)\;.
\]

The order $2$ element $S=\begin{bmatrix} \ \ 0 & 0&-1\\ \ \ 0&1&\ \ 0\\-1 & 0 &
\ \ 0 \end{bmatrix}\in\Ga$ acts by the map $(w_0,w)\mapsto (\frac{1}{w_0},
-\frac{w}{w_0})$ on $\hdc$. It thus fixes the point $x^\bullet=(1,0)$
and acts by $-\id$ on $T_{x^\bullet}\hdc$. In particular it maps
$v^\bullet$ to $-v^\bullet$. By equivariance, we thus have
\[
\forall\;g\in G,\quad\iota\,\wt\varphi(g)=\wt\varphi(gS)\;.
\]
The element $S$ centralises $M$ and normalises $\Phi^\RR\,$; more
precisely,
\[
\forall\;t\in \RR, \quad S\Phi^tS^{-1}=\Phi^{-t}\;.
\]
Since $S$ is the projective image of the matrix of the Hermitian form
$q=-z_0\overline{z_2} -z_2\overline{z_0}+|z_1|^2$, we have $^*g\,S\,g=S$
for every $g\in G$, hence
\[
\forall\;g\in G,\quad \;^*g^{-1}=S\,g\,S^{-1}\;.
\]
For all $x\in \Ga\bs G$ and $s\in\RR$, we again have $\Theta(x\Phi^s)
=\Theta(x)\Phi^{-s}$ and $\iota\circ \varphi=\varphi\circ \Theta$.

The (closed) horoball in $\hdc$ centered at $\infty$ whose boundary
$\partial \H_\infty$ contains $x^\bullet$ is
\[
\H_\infty=\{(w_0,w)\in \hdc:2\,\Re\; w_0 -|w|^2\geq 2\}\;.
\]
The Heisenberg group $H$ acts simply transitively on $\partial
\H_\infty$ and on $\normalpm \H_\infty$, which contains $\pm
v^\bullet$. Thus Equation \eqref{eq:relatnormalinoutwtphi} is still
satisfied. By for instance \cite[page 90]{ParPau17MA}, the stabiliser
$\Ga\!_{\H_\infty}$ in $\Ga$ of the horoball $\H_\infty$, as well as the
one of $\normalpm \H_\infty$, is equal to $\Ga_H$.  The
$\Ga$-equivariant families of horoballs
\[
\D^-=\D^+=(\ga\H_\infty)_{\ga\in\Ga}
\]
are again locally finite, since $\infty$ is a bounded parabolic fixed
point of $\Ga$.

For every $\ga\in\Ga$ having a representative (whose choice does not
change the following claims) in $\SU_q$ with first column
$\begin{pmatrix} a \\ \alpha \\ c \end{pmatrix}\in\M_{3,1}(\OOO_K)$,
we have $\ga\notin \Ga\!_{\H_\infty}$ if and only if $c\neq 0$ (see for
instance \cite[Eqs.~(42)]{ParPau10GT}) and then
\begin{enumerate}
\item[(i)] since $\infty=[1:0:0]$, the point at infinity $\ga
  \cdot\infty$ is equal to $\big(\frac{a}{c},\frac{\alpha}{c}\big)$;
\item[(ii)] since $H$ acts simply transitively on
  $\partial_\infty\hdc-\{\infty\}$, there exists a unique $r_\ga\in H$
  such that $r_\ga\cdot 0=\ga \cdot\infty$, and we have
  $r_\ga=n_-\big(\frac{a}{c},\frac{\alpha}{c}\big)$;
\item[(iii)] by \cite[Lem.~6.3]{ParPau10GT}, we have
  $d(\H_\infty,\ga\H_\infty)=\ln|c|=\frac{1}{2}\ln(\n(c))$.
\end{enumerate}
Therefore by \cite[Prop.~6.5 (2)]{ParPau10GT} with $\I=\OOO_K$, the
map $\ga\mapsto r_\ga$ induces, for all $t,t_0\in\RR$, a bijection
from $\{[\ga]\in\Ga\!_{\H_\infty}\bs\Ga/\Ga\!_{\H_\infty}:
0<d(\H_\infty,\ga\H_\infty)\leq t-t_0\}$ to $\F_{t-t_0}$.

Again using the simple transitivity of the action of $H$ on
$\normalpm\H_\infty$, we have $\Ga_H$-equivariant homeomorphisms
$\wt\psi_\pm: \normalpm\H_\infty\ra H$ which associates to
$v\in\normalpm\H_\infty$ the unique element $\wt\psi_\pm(v)\in H$
such that $\wt\psi_\pm(v)\cdot (\pm v^\bullet)=v$.

For every $\ga\in\Ga-\Ga\!_{\H_\infty}$, with $\rho_\ga$ the element of
$\normalout\H_\infty$ whose point at infinity is $\ga\cdot\infty$, the
maps $\wt\psi_\pm$ induce homeomorphisms $\psi_\pm:
\Ga_H\bs\normalpm\H_\infty \ra\Ga_H\bs H$ such that
\[
(\psi_+)_*(\Delta_{\Ga\rho_\ga})=\Delta_{\Ga_H r_\ga}\;.
\]
In the remainder of the proof of Corollary \ref{coro:HeisenbergFarey},
we use the same normalisation of the Patterson-Sullivan measures
$(\mu_x)_{x\in\hdc}$ as in \cite[\S 4]{ParPau17MA}. We denote by
$\delta_{x,y}=\Big\{\begin{array}{l} 1\;{\rm if}\;x=y\\0\;{\rm
  otherwise}\end{array}$ the Kronecker symbol.

\blemm\label{lem:calcskinHeisenberg} We have $\|\sigma^\mp_{\D^\pm}\|=
\frac{(1+2\,\delta_{D_K,-3})\,|D_K|}{4\,|\OOO_K^\times|}$.
\elemm

\dem By \cite[Lem.~12 (iv)]{ParPau17MA} with $n=2$, we have
$\|\sigma^\mp_{\D^\pm}\|=8\Vol(\Ga\!_{\H_\infty}\bs \H_\infty)$, where
$\Vol$ is the Riemannian volume. Denoting as in \cite[\S
  3]{ParPau17MA}, for every $s\in\RR$,
\[
\H_s=\{(w_0,w)\in \hdc:2\,\Re\; w_0 -|w|^2\geq s\}\;,
\]
we have $\H_\infty=\H_2$ and the horoballs $\H_s$ all have the same
stabiliser $\Ga\!_{\H_s}=\Ga\!_{\H_\infty}$ for $s\in\RR$. By the comment
following \cite[Eq.~(11)]{ParPau17MA}, we have $\Vol(\Ga\!_{\H_\infty}
\bs \H_\infty)=\frac{1}{4}\,\Vol(\Ga\!_{\H_1}\bs\H_1)$. The result then
follows from \cite[Lem.~16]{ParPau17MA} which says that
\[
\Vol(\Ga\!_{\H_1} \bs \H_1)=\frac{(1+2\delta_{D_K,-3})\,|D_K|}
{8\,|\OOO_K^\times|}\;.\quad\Box
\]

\medskip
Since we normalised $\mu_{\Ga_H\bs H}$ to be a probability measure, it
follows from Lemma \ref{lem:calcskinHeisenberg} that for $x\in\Ga_H\bs
H$,
\[
(\psi_+)_*(\sigma^+_{\D^-})=(\varphi^{-1})_*(\sigma^+_{\D^-})=
\frac{(1+2\,\delta_{D_K,-3})\,|D_K|}{4\,|\OOO_K^\times|}\;
\mu_{\Ga_H\bs H}\;.
\]

By \cite[Lem.~12 (iii)]{ParPau17MA} with $n=2$ and by the volume
formula of Holzapfel-Stover (see \cite[Lem.~17]{ParPau17MA} for the
appropriate normalisation of the volume form), we have
\[
\|m_{\rm BM}\|=\frac{\pi^2}{2}\;\Vol(M)=
\frac{\pi\,(1+2\,\delta_{D_K,-3})\,|D_K|^{5/2}\,\zeta_K(3)}{96\;\zeta(3)}\;.
\]

By \cite[Eq.~(21)]{ParPau17MA} and the comment following it, the index
of $H\cap \Ga$ in $\Ga_H$ is equal to $\frac{|\OOO_K^\times|}
{1+2\,\delta_{D_K,-3}}$. The map from $\Big\{(a,\alpha,c)\in\OOO_K
\times \OOO_K \times\OOO_K: \begin{array}{c}\langle a,\alpha,c\rangle
=\OOO_K\\ \tr(a\,\overline{c})=\n(\alpha),c\neq 0\end{array}\Big\}$
to $H$ defined by $(a,\alpha,c)\mapsto n_-(\frac{a}{c},
\frac{\alpha}{c})$ is $|\OOO_K^\times|$-to-$1$ onto its image.
Hence, using the (lifted linear) action of $n_-(w_0,w)\in H\cap \Ga$
on $(a,\alpha,c)\in\OOO_K \times \OOO_K \times\OOO_K$ defined by
\[
n_-(w_0,w)\cdot(a,\alpha,c)=
(a+\overline{w}\,\alpha +w_0\,c,\alpha+\omega \,c,c)\,,
\]
by \cite[Theo.~4]{ParPau17MA}, for every $t_0\in\RR$, we have, as
$t\ra+\infty$,
\begin{align*}
  \card\;\F_{t-t_0}&=\frac{1+2\,\delta_{D_K,-3}}{|\OOO_K^\times|^2}
  \times\\&\qquad
\card\Big(\,(H\cap \Ga)\bs\Bigg\{(a,\alpha,c)\in\OOO_K
\times \OOO_K \times\OOO_K: \begin{array}{c}\langle a,\alpha,c\rangle
  =\OOO_K\\ \tr(a\,\overline{c})=\n(\alpha)\\0<\n(c)\leq e^{2\,t-2\,t_0}
\end{array}\Bigg\}\,\Big)
\\
&\sim \frac{3\,(1+2\,\delta_{D_K,-3})\,\zeta(3)}
  {2\,\pi\,|\OOO_K^\times|^2\,\sqrt{|D_K|}\;\zeta_K(3)}\;e^{4\,t-4\,t_0}\;.
\end{align*}

Since $H\cap \Ga$ has finite index in $\Ga_H=\Ga\!_{\H_\infty}$ and acts
freely on $\partial \H_\infty$, there are only finitely many elliptic
elements in $\Ga$ up to conjugation by $\Ga\cap H$ whose fixed point
set contains $\infty=[1:0:0]$ as a point at infinity.  There are only
finitely many $\Ga\!_{\H_\infty}$-orbits of images of $\H_\infty$ by
$\Ga$ meeting $\H_\infty$. Hence there again exists a finite subset
$F$ of the set of double cosets $\Ga\!_{\H_\infty}\bs\Ga/
\Ga\!_{\H_\infty}$ such that for every $[\ga]\in\Ga\!_{\H_\infty}\bs\Ga/
\Ga\!_{\H_\infty}-F$, we have
\[
d(\H_\infty,\ga\H_\infty)>0 \quad {\rm and}\quad m_\ga=1\;.
\]
We have similarly, for all $\ga\in \Ga-\Ga\!_{\H_\infty}$ and $t\in\RR$,
\[
(\varphi^{-1})_*(\Delta_{\Ga\flow{t}\rho_\ga})=
\Delta_{\Ga r_\ga \Phi^{t}M}
\]
and by Lemma \ref{lem:calcskinHeisenberg}, for all $y\in \Ga_H\bs H$
and $s\in\RR$ with $s\geq t_0$,
\begin{align*}
d\big((\varphi^{-1})_*(\mu^{0+}_{\D^+,t_0})\big)
(\Theta(y\,\Phi^{-s})) &=\|\sigma^-_{\D^+}\|\;
d\mu_{\Ga_H\bs H}(y)\,e^{-4\,s}\,ds\\ &
=\frac{(1+2\,\delta_{D_K,-3})\,|D_K|}{4\,|\OOO_K^\times|}\;
d\mu_{\Ga_H\bs H}(y)\,e^{-4\,s}\,ds\;.
\end{align*}
The end of the proof of Corollary \ref{coro:HeisenbergFarey} proceeds
now as the one of Corollary \ref{coro:marklof}.
\cqfd

\subsection{Equidistribution of quaternionic Heisenberg Farey fractions
  at a given density}
\label{subsec:heisenbergFareyquat}

In this section, we denote by $\HH$ Hamilton's quaternion algebra over
$\RR$, with $x\mapsto \overline{x}$ its conjugation, $\n: x\mapsto
x\overline{x}$ its reduced norm, $\tr: x\mapsto x+\overline{x}$ its
reduced trace. Let $A$ be a definite ($A\otimes_\QQ\RR=\HH$)
quaternion algebra over $\QQ$, with discriminant $D_A$.  Let $\OOO$ be
a maximal order in $A$, with $\OOO^\times$ its finite group of
invertible elements. We denote by $_\OOO\langle a,\alpha,c\rangle$ the
left ideal of $\OOO$ generated by $a,\alpha,c\in \OOO$.  See
\cite{Vigneras80} for definitions.

Let $q$ be the nondegenerate quaternionic Hermitian form
 of Witt signature $(1,2)$ on the
right vector space $\HH^3$ over $\HH$ with coordinates
$(z_0,z_1,z_2)$ defined by
\[
q=-\tr(\overline{z_0}\,z_2) +\n(z_1)\;.
\]
With $\operatorname{U}_q=\{g\in \GL_3(\HH): q\circ g=q\}$, let
$G=\PU_q=\operatorname{U}_q/\{\pm\id\}$ be the projective unitary group
of $q$.  Let $\Ga$ be the image of $\operatorname{U}_q\cap\GL_3(\OOO)$
in $G$, which is a (nonuniform) arithmetic lattice in $G$.

Denoting by $\begin{bmatrix} a & \overline{\ga} & b \\
  \alpha & A & \beta\\c & \overline{\delta} & d \end{bmatrix}$
the image in $G$ of $\begin{pmatrix} a & \overline{\ga} & b \\
  \alpha & A & \beta\\c & \overline{\delta} & d \end{pmatrix}\in
\operatorname{U}_q$, let
\[
H=\Big\{n_-(w_0,w)=\begin{bmatrix} 1 & \overline{w} & w_0\\0&1&w
\\0&0&1\end{bmatrix}:w_0,w\in\HH, \;\tr(w_0)=\n(w)\Big\}\,,
\]
\[
\Phi^\RR=\Big\{
\Phi^t=\begin{bmatrix} \;e^{-t} & 0& 0\\0&1&0\\0&0&e^{t}\end{bmatrix}
:t\in\RR\Big\}\quad {\rm and}\quad
\]
\[M=\Big\{m(u,U)=
\begin{bmatrix} \;u & 0& 0\\0&U&0\\0&0&u\end{bmatrix}
:u,U\in\HH,\;\n(u)=\n(U)=1\Big\}\,.
\]
Since $\RR$ is central in $\HH$, the subgroup $M$ is the compact
factor of the centraliser in $G$ of the standard Cartan subgroup
$\Phi^\RR$ of $G$, and the subgroup $M\Phi^\RR$ normalises the {\em
  quaternionic Heisenberg group} $H$, since
\[
m(u,U)\,n_-(w_0,w)\,m(u,U)^{-1}=
n_-(u\,w_0\,\overline{u},U\,w\,\overline{u})\;.
\]
Since $\OOO$ is invariant under conjugation in $\HH$, the groups $\Ga$
and $M$ are invariant under the standard Cartan involution $$g\mapsto
\;^*g^{-1}\,,$$ where $^*g$ is the image in $G$ of the
transpose-conjugate matrix of any matrix in $\operatorname{U}_q$
representing $g$.

Let $\Ga_H=N_G(H)\cap\Ga= (MH)\cap \Ga= \Big\{\begin{bmatrix} 
\;u &u\,\overline{v}& u\,v_0\\
0&U&U\,v\\
0&0&u\end{bmatrix}:
\begin{array}{c} u,U\in\OOO^\times, \;v,v_0\in\OOO\\\tr(v_0)=\n(v)
\end{array}\Big\}$, which acts properly discontinuously on the left
on $H$ by (noting the lack of commutativity)
\begin{equation}\label{eq:actionGasubHsurHquat}
  \begin{bmatrix} \;u & u\,\overline{v}& u\,v_0\\0&U&U\,v\\0&0&u
  \end{bmatrix}\cdot\begin{bmatrix} \;1 &
\overline{w}& w_0\\0&1&w\\0&0&1\end{bmatrix}
=\begin{bmatrix} 
\;1 &u(\overline{w}+\overline v)\,\overline U&
u(v_0+w_0+\overline{v}\,w)\overline u\,\\
0&1&U\,(w+v)\overline u\\
0&0&1\end{bmatrix}\,.
\end{equation}
The inclusion map $H\ra G$ again induces an identification between the
quotient $\Ga_H\bs H$ and the image of $H$ in $\Ga\bs G/M$.  We again
endow $\Ga_H\bs H$ with the induced measure $\mu_{\Ga_H\bs H}$ of a
Haar measure on $H$, normalised to be a probability measure, that we
also see as a probability measure on $\Ga\bs G/M$ (with support
$\Ga_H\bs H$).

For every $t\in\RR$, we consider the subset $\F_{t}$ of $\Ga_H\bs H$
consisting of the {\em quaternonic Heisenberg Farey fractions of
  height at most $e^{t}$}, defined by
\[
\F_{t}=\Ga_H\bs\Big\{n_-(a\,c^{-1},\alpha\,c^{-1}):
\begin{array}{c}a,\alpha,c\in\OOO,\;\;_\OOO\langle a,\alpha,c\rangle
  =\OOO,\\ \tr(\overline{a}\,c\,)=\n(\alpha),\end{array}\quad
0<\n(c)\leq e^{2\,t}\Big\}\,.
\]
Note that the above set of elements $n_-(a\,c^{-1},\alpha\,c^{-1})$ is
indeed invariant under $\Ga_H$, by Equation
\eqref{eq:actionGasubHsurHquat}.  Let $\Theta:\Ga\bs G/M\ra \Ga\bs
G/M$ be the Cartan involutive homeomorphism defined by $\Ga gM\mapsto
\Ga{\;}^tg^{-1}M$.

\bcoro\label{coro:HeisenbergFareyquat} For every $t_0\in\RR$, for the
narrow convergence of probability measures on $(\Ga_H\bs H)
\times(\Ga\bs G/M)$, we have
\begin{align*}
&\lim_{t\to+\infty}\;\frac 1{\card\,\F_{t-t_0}}
\sum_{r\in\F_{t-t_0}}\Delta_r\otimes\Delta_{\Ga \,r\,\Phi^tM}\\=\;&
10\;e^{10\,t_0} \int_{s=t_0}^{+\infty} (\mu_{\Ga_H\bs H})\otimes
\big(\Theta_*\,(\Phi^{-s})_*\,\mu_{\Ga_H\bs H}\big)\,e^{-10\,s}\,ds\,.
\end{align*}
\ecoro

As a remark similar to the remark at the end of Section
\ref{subsec:retrouvonsMarklof}, one could obtain an error term under
an additional smoothness assumption, and a joint partial
equidistribution result of quaternionic Heisenberg Farey points
$n_-\big(a\,c^{-1}, \alpha\,c^{-1})$ modulo $\Ga_H$ with their
denominators $c$ congruent to $0$ modulo any fixed element $N$ in
$\OOO-\{0\}$.

\medskip
\dem We mostly indicate the differences with the proof of Corollary
\ref{coro:HeisenbergFarey}. We refer to \cite{Mostow73,KimPar03,
  Philippe16} as well as \cite[\S 3]{ParPau22MPCPS} for background on
quaternionic hyperbolic geometry. We follow the conventions of this
last reference concerning the normalisation of the sectional curvature
and the choice of the quaternionic Hermitian form with Witt signature
$(1,2)$.

We now consider $X=\hdh$ the Siegel domain model of the quaternionic
hyperbolic plane, that is, the quaternionic manifold
\[
\big\{(w_0,w)\in\HH^2: \tr (w_0) -\n(w)>0\big\}\,,
\]
endowed with the Riemannian metric
\begin{equation}\label{eq:metriemSiegeldomquat}
ds^2_{\,\hdh}=\frac{1}{(\tr w_0 -\n(w))^2}
\big(\n(dw_0-\overline{dw}\, w)+
(\tr(w_0) -\n(w))\,\n(dw)\big)\,.
\end{equation}
This metric is again normalised so that its sectional curvatures are
in $[-4,-1]$. The boundary at infinity of $\hdh$ is
\[
\partial_\infty\hdh=\big\{(w_0,w)\in
\HH^2 \;:\tr( w_0) -\n(w)=0\big\}\cup\{\infty\}\,.
\]
Using right-homogeneous coordinates, we identify $\hdh\cup
\partial_\infty \hdh$ with its image in the right projective plane
$\PP^2_r(\HH)$ over $\HH$ by the map $(w_0,w)\mapsto [w_0:w:1]$ and
$\infty \mapsto [1:0:0]$. We denote by $\cdot$ the left projective
action of $G$ on $\hdh\cup\partial_\infty \hdh$, as well as its
derived action on $T^1\hdh$.

The critical exponent of the (nonuniform arithmetic) lattice $\Ga$ of
$G$ is now (see for instance \cite[Theo.~4.4 (i)]{CorIoz99})
\[
\delta_\Ga=10\;.
\]

We again fix $v^\bullet=((1,0),(-2,0))\in T^1\hdh$, which is indeed a
unit tangent vector with footpoint $x^\bullet =(1,0)$ by Equation
\eqref{eq:metriemSiegeldomquat}. The stabiliser of $v^\bullet$ in $G$
is again equal to $M$ and is hence centralised by $\Phi^\RR$. The
$G$-equivariant orbital map $\wt\varphi: g\mapsto g\cdot v^\bullet$
now defines a homeomorphism $\varphi: \Ga\bs G/M \ra \Ga\bs T^1\hdh$.

For every $t\in\RR$, the element $\Phi^t$ acts on $\hdh$ by the map
$(w_0,w)\mapsto (e^{-2t}w_0,e^{-t}w)$. The geodesic line $\ell$ in
$\hdh$ such that $\ell(0)=x^\bullet$ and $\ell'(0)= v^\bullet$ is
$t\mapsto (e^{-2\,t},0)$. Hence, as in the complex case (see the proof
of Corollary \ref{coro:HeisenbergFarey}), we have
\[
\forall\;t\in\RR,\;\forall\;g\in G,\quad \flow{t}\wt\varphi(g)=
\wt\varphi(g\Phi^t)\;.
\]

The order $2$ element $S=\begin{bmatrix} \ \ 0 & 0&-1\\\ \ 0&1&\ \ 0\\-1 & 0 &
\ \ 0 \end{bmatrix}$ still belongs to $\Ga$, it centralises $M$ and
normalises $\Phi^\RR$, and it acts by the map $(w_0,w)\mapsto
(w_0^{-1},-w\,w_0^{-1})$ on $\hdh$. Since $S$ is the projective image
of the matrix of the quaternionic Hermitian form
$q=-\tr(\overline{z_0}\,z_2) +\n(z_1)$, we have
\[
\forall\;g\in G,\quad \;^*g^{-1}=S\,g\,S^{-1}\;.
\]
As in the complex case, for all $g\in G$, $t\in\RR$ and $x\in \Ga\bs
G/M$, we have
\[
\iota\,\wt\varphi(g)=\wt\varphi(gS),\quad S\Phi^tS^{-1}=\Phi^{-t},
\quad \iota\circ \varphi=\varphi\circ\Theta
\quad{\rm and}\quad \Theta(x\Phi^t)=\Theta(x)\Phi^{-t}\;.
\]

The (closed) horoball in $\hdh$ centered at $\infty$ whose boundary
$\partial \H_\infty$ contains $x^\bullet$ is
\[
\H_\infty=\{(w_0,w)\in \hdh:\tr( w_0) -\n(w)\geq 2\}\;.
\]
The quaternionic Heisenberg group $H$ again acts simply transitively on
$\partial \H_\infty$, and on $\normalpm \H_\infty$ which contains $\pm
v^\bullet$. Thus Equation \eqref{eq:relatnormalinoutwtphi} is still
satisfied. By for instance the end of \S 3 in \cite{ParPau22MPCPS},
the stabiliser $\Ga\!_{\H_\infty}$ in $\Ga$ of the horoball $\H_\infty$,
as well as the one of $\normalpm \H_\infty$, is equal to $\Ga_H$.  The
$\Ga$-equivariant families of horoballs
\[
\D^+=\D^-=(\ga\H_\infty)_{\ga\in\Ga}
\]
are again locally finite, since $\infty$ is again a bounded parabolic
fixed point of $\Ga$.

For every $\ga\in\Ga$ having a representative in $U_q$ with first
column $\begin{pmatrix} a \\ \alpha\\ c \end{pmatrix}\in
\M_{3,1}(\OOO)$, we have $\ga\notin \Ga\!_{\H_\infty}$ if and only if
$c\neq 0$ (see for instance \cite{KimPar03},
\cite[Eqs.~(3$\cdot$3)]{ParPau22MPCPS}) and then
\begin{enumerate}
\item[(i)] since $\infty=[1:0:0]$, the point at infinity $\ga
  \cdot\infty$ is equal to $(a\,c^{-1},\alpha\,c^{-1})$;
\item[(ii)] since $H$ acts simply transitively on
  $\partial_\infty\hdc-\{\infty\}$, there exists a unique $r_\ga\in H$
  such that $r_\ga\cdot 0=\ga \cdot\infty$, and we have $r_\ga=
  n_-(a\,c^{-1},\alpha\,c^{-1})$;
\item[(iii)] with $\H_s =\{(w_0,w)\in\hdh: \tr(w_0)-\n(w)=s\}$ for $s>0$, by
  \cite[Lem.~6$\cdot$5]{ParPau22MPCPS} where we take $s=2$ so that
  $\H_s=\H_\infty$, we have $d(\H_\infty, \ga\H_\infty)=
  \frac{1}{2}\ln(\n(c))$.
\end{enumerate}
Therefore by \cite[Prop.~4$\cdot$2 (ii)]{ParPau22MPCPS} with $\mmm
=\OOO$, the map $\ga\mapsto r_\ga$ induces, for all $t,t_0 \in\RR$, a
bijection from $\{[\ga]\in\Ga\!_{\H_\infty}\bs\Ga/ \Ga\!_{\H_\infty}:
0<d(\H_\infty,\ga\H_\infty)\leq t-t_0\}$ to $\F_{t-t_0}$.

As in the complex case, we have homeomorphisms $\psi_\pm:
\Ga_H\bs\normalpm\H_\infty \ra\Ga_H\bs H$ such that
\[
(\psi_+)_*(\Delta_{\Ga\rho_\ga})=\Delta_{\Ga_H r_\ga}\;.
\]
In the remainder of the proof of Corollary \ref{coro:HeisenbergFareyquat},
we use the same normalisation of the Patterson-Sullivan measures
$(\mu_x)_{x\in\hdh}$ as in \cite[\S 7]{ParPau22MPCPS}. 

\blemm\label{lem:calcskinHeisenbergquat} We have $\displaystyle
\|\sigma^\mp_{\D^\pm}\|=\frac{D_A^2}{64\;|\OOO^\times|^2}$.
\elemm

\dem By \cite[Lem.~7$\cdot$2 (iv)]{ParPau22MPCPS} with $n=2$, we have
$\|\sigma^\mp_{\D^\pm}\|=80\Vol(\Ga\!_{\H_\infty}\bs \H_\infty)$, where
$\Vol$ is the Riemannian volume.  By
\cite[Lem.~7$\cdot$1]{ParPau22MPCPS} and the arguments in its proofs,
and by Equation (8.4) of loc.~cit. for the last equality, we have
\begin{align*}
\Vol(\Ga\!_{\H_\infty} \bs \H_\infty)&=\frac{1}{10}\Vol(\Ga\!_{\H_\infty} \bs
\partial\H_\infty)=\frac{1}{10}\,\frac{1}{2^5}\,\Vol(\Ga\!_{\H_1} \bs
\partial\H_1))=\frac{1}{2^5}\,\Vol(\Ga\!_{\H_1} \bs\H_1)\\&=
\frac{1}{2^5}\,\frac{D_A^2}{160\;|\OOO^\times|^2}\;.
\end{align*}
The result follows. \cqfd

\medskip
Since we normalised $\mu_{\Ga_H\bs H}$ to be a probability measure, it
follows from Lemma \ref{lem:calcskinHeisenbergquat} that for $x\in\Ga_H\bs
H$,
\[
(\psi_+)_*(\sigma^+_{\D^-})=(\varphi^{-1})_*(\sigma^+_{\D^-})=
\frac{D_A^2}{64\;|\OOO^\times|^2}\;\mu_{\Ga_H\bs H}\;.
\]

Let $m_A=24$ if $D_A$ is even, and $m_A=1$ otherwise. By respectively
Lemma 7$\cdot$2 (iii) with $n=2$ and Theorem 1$\cdot$4 in
\cite{ParPau22MPCPS}, we have, with $p$ ranging over primes,
\[
\|m_{\rm BM}\|=\frac{\pi^4}{48}\;\Vol(M)=
\frac{\pi^8\,m_A}{2^{17}\cdot3^6\cdot5^2\cdot7}
\;\prod_{p\,\mid\,D_A}(p-1)(p^2+1)(p^3-1)\;.
\]

By the definition of $\Ga_H$, the index of
$H\cap \Ga$ in $\Ga_H$ is now equal to $\frac{|\OOO^\times|^2}{2}$.
The map from $\Big\{(a,\alpha,c) \in\OOO \times \OOO \times\OOO
: \begin{array}{c} _\OOO\langle a,\alpha,c\rangle =\OOO
  \\ \tr(\overline{a}\,c) =\n(\alpha),c\neq 0\end{array}\Big\}$ to $H$
  given by $(a,\alpha,c)\mapsto n_-(a\,c^{-1}, \alpha\,c^{-1})$ is
$|\OOO^\times|$-to-$1$ onto its image.  Hence, using the (lifted
linear) action of $n_-(w_0,w)\in H\cap \Ga$ on $(a,\alpha,c)\in\OOO
\times \OOO\times\OOO$ defined by
\[
n_-(w_0,w)\cdot(a,\alpha,c)=
(a+\overline{w}\,\alpha +w_0\,c,\alpha+\omega \,c,c)\,,
\]
by \cite[Theo.~1$\cdot$1]{ParPau22MPCPS}, for every $t_0\in\RR$, we
have, as $t\ra+\infty$,
\begin{align*}
  \card\;\F_{t-t_0}&=\frac{2}{|\OOO^\times|^3}\;
\card\Big((H\cap \Ga)\bs\Big\{(a,\alpha,c)\in\OOO
\times \OOO \times\OOO: \begin{array}{c}_\OOO\langle a,\alpha,c\rangle
  =\OOO\\ \tr(\overline{a}\,c)=\n(\alpha)\\0<\n(c)\leq e^{2\,t-2\,t_0}
\end{array}\Big\}\Big) \\&\sim \frac{2^4\cdot3^6\cdot5\cdot7\;D_A^4}
{\pi^8 \;m_A\;|\OOO^\times|^4\;\prod_{p\,\mid\, D_A}(p-1)(p^2+1)(p^3-1)}
\;e^{10\,t-10\,t_0}\;.
\end{align*}

As in the complex case, there exists a finite subset $F$ of
$\Ga\!_{\H_\infty}\bs\Ga/ \Ga\!_{\H_\infty}$ such that for every
$[\ga]\in\Ga\!_{\H_\infty}\bs\Ga/ \Ga\!_{\H_\infty}-F$, we have
\[
d(\H_\infty,\ga\H_\infty)>0, \quad m_\ga=1, \quad 
(\varphi^{-1})_*(\Delta_{\Ga\flow{t}\rho_\ga})=
\Delta_{\Ga r_\ga \Phi^{t}M}
\]
and by Lemma \ref{lem:calcskinHeisenbergquat}, for all $y\in \Ga_H\bs
H$ and $s\in\RR$ with $s\geq t_0$,
\begin{align*}
d\big((\varphi^{-1})_*(\mu^{0+}_{\D^+,t_0})\big)
(\Theta(y\,\Phi^{-s})) &=\|\sigma^-_{\D^+}\|\;
d\mu_{\Ga_H\bs H}(y)\,e^{-10\,s}\,ds\\ &
=\frac{D_A^2}{64\;|\OOO^\times|^2}\;
d\mu_{\Ga_H\bs H}(y)\,e^{-10\,s}\,ds\;.
\end{align*}
The end of the proof of Corollary \ref{coro:HeisenbergFareyquat} proceeds
now as the one of Corollary \ref{coro:HeisenbergFarey}.
\cqfd

\subsection{Equidistribution of nonarchimedian Farey fractions at a
  given density}
\label{subsec:nonarchFarey}

In this section, we give an arithmetic application of Theorem
\ref{theo:quotient} (3), proving a joint partial equidistribution
result of nonarchimedean arithmetic points with given density on an
expanding horosphere in the quotient of a regular tree by a nonuniform
arithmetic lattice.

Let $K$ be a (global) function field of genus $\ggg$ over a finite
field $\FF_q$ of order a positive prime power $q$, let $v$ be a
(normalised discrete) valuation of $K$, let $K_v$ be the associated
completion of $K$, let $\OOO_v=\{x\in K_v:v(x)\geq 0\}$ be its
valuation ring, let $\pi_v\in K$ with $v(\pi_v)=1$ be a uniformiser of
$v$, let $q_v$ be the order of the residual field
$\OOO_v/\pi_v\OOO_v$, let $|\cdot|_v=q_v^{-v(\,\cdot\,)}$ be the
absolute value associated with $v$, and let $R_v$ be the affine
function ring associated with $v$. For all these notions and
complements, we refer to \cite{Goss96, Rosen02}, as well as to
\cite[\S 14.2]{BroParPau19} whose notation we will follow. The
simplest example, used in Corollary \ref{coro:nonarchFareyintro}, is
given by $K=\FF_q(Y)$ the field of rational fractions over $\FF_q$
with one indeterminate $Y$, $\ggg=0$, $v=v_\infty:\frac{P}{Q}\mapsto
\deg Q-\deg P$ for every $P,Q\in \FF _q[Y]$ the valuation at infinity,
$K_v=\FF_q((Y^{-1}))$, $\OOO_v= \FF_q[[Y^{-1}]]$ the local ring of
formal power series in $Y^{-1}$, $\pi_v=Y^{-1}$, $q_v=q$, and
$R_v=\FF_q[Y]$.

Let $G$ be the locally compact group $\PGL_2(K_v)=\GL_2(K_v)/
(K_v^\times\id)$. We denote by $\begin{bmatrix} a & b \\ c &
  d \end{bmatrix}$ the image in $G$ of $\begin{pmatrix} a & b \\ c &
  d \end{pmatrix}\in \GL_2(K_v)$. Let $\Ga=\PGL_2(R_v)$ be the {\it
  Nagao lattice} in $G$ (see for instance \cite{Weil70}). We consider
the subgroups of $G$ defined by
\[
H=\Big\{n_-(r)=\begin{bmatrix} 1 & r\\0&1\end{bmatrix}:r\in K_v\Big\},
\quad  \Phi^\ZZ=\Big\{\Phi^n=
\begin{bmatrix} \;1 & 0\\0&\pi_v^{-n}\end{bmatrix}:n\in \ZZ\Big\}\;,
\]
and $M=\Big\{\begin{bmatrix} \;1 & 0\\0&u\end{bmatrix}:u\in K_v,
|u|_v=1\Big\}$. Note that $M$ centralises the standard Cartan subgroup
$\Phi^\ZZ$, that the diagonal subgroup $M\Phi^\ZZ$ normalises $H$, and
that both $\Ga$ and $M$ are invariant under the standard Cartan
involution $g\mapsto \;^tg^{-1}$.

Let
\[
\Ga_H=N_G(H)\cap\Ga= (HM)\cap \Ga= \Big\{\begin{bmatrix} 1 &
b\\0&d\end{bmatrix}: d\in R_v^\times, b\in R_v\Big\}\,,
\]
which acts properly discontinuously on the left on $H$ by
$$\begin{bmatrix} 1 &b\\0&d\end{bmatrix}\cdot
\begin{bmatrix} 1 &r\\0&1\end{bmatrix}=
\begin{bmatrix} 1 &\frac{r+b}{d}\\0&1\end{bmatrix}\,.
$$
The inclusion map $H\ra G$ again induces an identification between the
quotient $\Ga_H\bs H$ and the image of $H$ in $\Ga\bs G/M$.  We again
endow $\Ga_H\bs H$ with the induced measure $\mu_{\Ga_H\bs H}$ of a
Haar measure on $H$, normalised to be a probability measure, that we
also see as a probability measure on $\Ga\bs G/M$ (with support
$\Ga_H\bs H$).

For every $n\in\ZZ$, we consider the subset $\F_{n}$ of $\Ga_H\bs H$
consisting of the {\em Farey fractions of
  height at most $q_v^n$ with respect to $v$}, defined by
\[
\F_{n}=\Ga_H\bs\Big\{n_-\big(\frac{a}{c}\big):
\begin{array}{c}a,c\in R_v,\;\; aR_v+cR_v=R_v\\
c\neq 0, \;\;v(c)\geq -n\end{array}\Big\}\,.
\]
Let $\Theta:\Ga\bs G/M\ra \Ga\bs G/M$ be the Cartan involutive
homeomorphism defined by $\Ga gM\mapsto \Ga{\;}^tg^{-1}M$.

\bcoro\label{coro:nonarchFarey} For every $n_0\in\ZZ$, for the narrow
convergence of probability measures on $(\Ga_H\bs H) \times(\Ga\bs
G/M)$, we have
\begin{align*}
\lim_{n\to+\infty}\;&\frac 1{\card\,\F_{n-n_0}}
\sum_{r\in\F_{n-n_0}}\Delta_r\otimes\Delta_{\Ga \,r\,\Phi^{2n}M}\\
&=(q^2_v-1)\;q_v^{2n_0-2} \sum_{m=n_0}^{+\infty}
(\mu_{\Ga_H\bs H})\otimes
\big(\Theta_*\,(\Phi^{-2m})_*\,\mu_{\Ga_H\bs H}\big)\,q_v^{-\,2m}\,.
\end{align*}
\ecoro

Corollary \ref{coro:nonarchFareyintro} follows by considering the
particular valued function field $(\FF_q(Y),v_\infty)$ indicated
above. As a remark similar to the remark at the end of Section
\ref{subsec:retrouvonsMarklof}, one could obtain an error term under
an additional locally constant regularity assumption, and a joint
partial equidistribution result of nonarchimedian Farey points
$n_-\big(\frac{a}{c}\big)$ modulo $\Ga_H$ with their denominators $c$
congruent to $0$ modulo any fixed element $N$ in $R_v-\{0\}$.

\medskip
\dem We mostly indicate the differences with the proof of Corollary
\ref{coro:HeisenbergFarey}. We refer to \cite{Tits79, Serre83} for
background on Bruhat-Tits trees, as well as to \cite[\S 15.1 and \S
  15.2]{BroParPau19} whose notation we will follow.

We now consider $\XX=\XX_v$ the Bruhat-Tits tree of $(\PGL_2,K_v)$,
which is a regular tree of degree $q_v+1$ endowed with a vertex
transitive action of $G$. Note that $\Ga$ acts without inversion on
$\XX_v$ by \cite[II.1.3]{Serre83}. The set of vertices of $\XX_v$ is
the set of homothety classes $[\Lambda]$ under $K_v^\times$ of
$\OOO_v$-lattices $\Lambda$ in the plane $K_v\times K_v$, and
$g[\Lambda]=[g\Lambda]$ for every $g\in G$.  We identify the boundary
at infinity $\partial_\infty \XX_v$ of (the geometric realisation of)
$\XX_v$ and the projective line $\PP_1(K_v)=K_v\cup\{\infty\}$ by the
unique homeomorphism such that the (continuous) extension to
$\partial_\infty \XX_v$ of the isometric action of $G$ on $\XX_v$ is
the projective action of $G$ on $\PP_1(K_v)$, that is, the action of
$G$ by homographies on $K_v\cup\{\infty\}$. We denote by $\cdot$ the
action of $G$ by homographies on $K_v\cup\{\infty\}$, as well as the
action of $G$ on the space $\G\XX_v$ of (discrete) geodesic lines in
$\XX_v$.

The critical exponent of the (nonuniform arithmetic) lattice $\Ga$ of
$G$ is now (see for instance \cite[Eq.~(15.8)]{BroParPau19})
\begin{equation}\label{eq:expocritnonarch}
\delta_\Ga=\ln q_v\;.
\end{equation}

The standard basepoint $x^\bullet$ of $\XX_v$ is the homothety class
$[\OOO_v\times\OOO_v]$ of the standard $\OOO_v$-lattice $\OOO_v\times
\OOO_v$ in $K_v\times K_v$. We consider the geodesic line $v^\bullet
\in\G\XX_v$ with $v^\bullet(0)=x^\bullet$, $v^\bullet(-\infty) =\infty
\in\PP_1(K_v)$ and $v^\bullet(+\infty)=0\in\PP_1(K_v)$. The stabiliser
of $v^\bullet$ in $G$ is again equal to $M$. The $G$-equivariant
orbital map $\wt\varphi: g\mapsto g\cdot v^\bullet$ now defines a
homeomorphism $\varphi: \Ga\bs G/M \ra \Ga\bs \G\XX_v$.

Since $v^\bullet(n)=[\OOO_v\times\pi_v^{-n}\OOO_v]$ for every
$n\in\ZZ$ (see for instance \cite[top of page 310]{BroParPau19}) and
by equivariance, we have (see also \cite[Eq.~(15.4)]{BroParPau19})
\[
\forall\; n\in\ZZ,\;\forall\;g\in G, \quad
\flow{n}\wt\varphi(g)=\wt\varphi(g\,\Phi^n)\;.
\]

The order $2$ element $S=\begin{bmatrix} 0 & -1\\ 1 & \ \ 0 \end{bmatrix}$
still belongs to $\Ga$, it normalises $M$ and $\Phi^\RR$, more
precisely $S\,\Phi^nS^{-1}=\Phi^{-n}$ for every $n\in\ZZ$. By
equivariance, the antipodal map $\iota$ satisfies $\iota\,
\wt\varphi(g) =\wt\varphi(gS)$ for every $g\in G$.  Since the
computation is independent of the ground field, we have $^tg^{-1} =
S\,gS^{-1}$ for every $g\in G$. Hence $\iota\circ\varphi= \varphi
\circ\Theta$ and $\Theta(x\Phi^n)=\Theta(x)\Phi^{-n}$ for all $x\in
\Ga\bs G/M$ and $n\in\ZZ$.

The group $H$ fixes the point at infinity $\infty$, preserves the
horoball $\H_\infty$ in $\XX_v$ centered at $\infty$ whose boundary
contains $x^\bullet$, and acts simply transitively on $\partial_\infty
\XX_v-\{\infty\}=K_v$, hence on $\normalpm\H_\infty$.  Note that
$\normalout\H_\infty$ contains the geodesic ray $v^\bullet{}_{\mid\,
[0,+\infty[}$ and that $\normalin\H_\infty$ contains $(\iota v^\bullet)
{}_{\mid\, ]-\infty,0]}$. In particular, we have $\normalout\H_\infty
=\{\ell_{\mid[0,+\infty[}:\ell\in W^-(v^\bullet)\}$.

Note that defining $V_{\rm even} \XX_v$, $\gengeod_{\rm even} \XX_v$
and $\G_{\rm even} \XX_v$ for the above basepoint $x^\bullet$ as just
before the statement of Theorem \ref{theo:Marklof6discrete}, we have
$\normalpm \H_\infty\subset \gengeod_{\rm even} \XX_v$, since any two
points of the horosphere $\partial \H_\infty$ are at even distance one
from the other. Furthermore, $\Ga$ preserves $V_{\rm even}
\XX_v$. Indeed, note that in a simplicial tree, if two of the
distances between three points are even, so is the third one. The
result then follows from \cite[II.1.2, Cor.]{Serre83}, which proves
that the distance $d(x^\bullet, \ga x^\bullet)$ is even for every
$\ga\in \GL_2(R_v)$, since $v(\det \ga) = 0$.

Each geodesic ray $w\in\normalin\H_\infty$ can be extended to a unique
element $\wh w\in\G\XX_v$ such that $\wh w(+\infty)$ is the point at
infinity of $\H_\infty$. This element belongs to $\G_{\rm even}\XX_v$,
is equal to $(N^+_{\iota \,v^\bullet})^{-1}(w)$ with the notation
$N^+_{\;\cdot}$ of Section \ref{sec:background}, and we define
$\widehat{\normalin\H_\infty}=\{\wh w:w\in \normalin\H_\infty\}$.
With $\Phi^{\geq n_0}=\{\Phi^n:n\geq n_0\}$, we have
$$
W_{n_0}^{0+}(\iota\, v^\bullet)=
\bigcup_{n\geq n_0}\flow{n}\;\widehat{\normalin\H_\infty}=
\bigcup_{n\geq n_0}\flow{n}H\,\iota\,v^\bullet=
\wt\varphi\big(H(\Phi^{\geq n_0})^{-1}S\big)\;.
$$

The subgroup $\Ga_H$ is again equal to the stabiliser
$\Ga\!_{\H_\infty}$ of the horoball $\H_\infty$ in $\Ga$, and $\infty$
is again a bounded parabolic fixed point of $\Ga$. We again consider
the locally finite $\Ga$-equivariant families of horoballs
\[
\D^+=\D^-=(\ga\H_\infty)_{\ga\in\Ga}\;.
\]
Note that the support of the skinning measure $\sigma^+_{\D^-}$ is
contained in $\Ga\bs\gengeod_{\rm even} \XX_v$, hence $\sigma^+_{\D^-}
\,_{\mid \;\Ga\bs\gengeod_{\rm even}\XX_v}=\sigma^+_{\D^-}$.

By \cite[Prop.~6.1]{Paulin02} when $K=\FF_q(Y)$ and $v=v_\infty$, and
by \cite[Lem.~15.1]{BroParPau19} in general, for every $\ga=
\begin{bmatrix} a & b\\c&d\end{bmatrix}\in\Ga$ with $c\neq 0$, we
have
\[
d(\H_\infty, \ga\H_\infty)=-2\,v(c)=2\,\log_{q_v}|c|_v\;.
\]
In particular, the distances $d(\H_\infty, \ga\H_\infty)$ for
$\ga\in\Ga$ are even and the endpoints of the common perpendiculars
between elements of $\D^-$ and $\D^+$ belong to $V_{\rm even} \XX_v$.
The map $\ga= \begin{bmatrix} a & b\\c&d\end{bmatrix}\mapsto
  n_-\big(\frac{a}{c}\big)$ now induces, for every $n\in\ZZ$, a
  bijection from
\[
\big\{[\ga]\in\Ga\!_{\H_\infty}\bs\Ga /\Ga\!_{\H_\infty}:
0<d(\H_\infty, \ga\H_\infty)\leq 2n\big\}
\]
to $\F_n$.  Denoting by $\rho_\ga$ the element of $\normalout
\H_\infty$ whose point at infinity is $\ga\cdot\infty=\frac{a}{c}$,
the map $\wt\psi_+: \normalout \H_\infty\ra H$ defined by $w\mapsto
n_-(w(+\infty))$ now induces an homeomorphism $\psi_+: \Ga_H\bs
\normalout\H_\infty \ra\Ga_H\bs H$, such that
\[
(\psi_+)_*(\Delta_{\Ga_H\rho_\ga})=\Delta_{\Ga_H n_-(\ga\cdot\infty)}\;.
\]

In the remainder of the proof of Corollary
\ref{coro:nonarchFarey}, we use the same normalisation of the
Patterson-Sullivan measures $(\mu_x)_{x\in V\XX_v}$ as in \cite[\S
  15.3]{BroParPau19}. Since we normalised $\mu_{\Ga_H\bs H}$ to be a
probability measure, it follows from \cite[Prop.~15.3
  (2)]{BroParPau19} that, for $x\in\Ga_H\bs H$,
\begin{equation}\label{eq:psipmstarnonarch}
  (\psi_+)_*(\sigma^+_{\D^-}) =(\varphi^{-1})_*(\sigma^+_{\D^-})
  =\frac{q^{\mathfrak g-1}}{q-1}\;\mu_{\Ga_H\bs H}\;.
\end{equation}

With $\zeta_K$ the Dedekind zeta function of $K$ (see for instance
\cite[\S 7.8]{Goss96} or \cite[\S 5]{Rosen02}), by \cite[Prop.~15.3
  (1)]{BroParPau19}, we have
\[
\|m_{\rm BM}\|=2\,\zeta_K(-1)\;\frac{q_v+1}{q_v}\;.
\]

By \cite[Eq.~(14.3)]{BroParPau19}, the subgroup $H\cap \Ga=n_-(R_v)$
has index $|R_v^\times|=q-1$ in $\Ga_H$.  The map from $\big\{(x,y)
\in R_v\times R_v : x R_v+yR_v=R_v, y\neq 0\big\}$ to $H$ given by
$(x,y) \mapsto n_-\big(\frac{x}{y}\big)$ is $|R_v^\times|$-to-$1$ onto
its image.  Hence, using the action by shears of $R_v $ on $R_v\times
R_v$ defined by $z\cdot(x,y)=(x+zy,y)$, by
\cite[Coro.~16.2]{BroParPau19} with $G=\GL_2(R_v)$ and
$(x_0,y_0)=(1,0)$ so that $m_{v,x_0,y_0}=q-1$ by
\cite[Eq.~(16.1)]{BroParPau19} with the notation of loc.~cit., for
every $n_0\in\ZZ$, as $n\ra+\infty$, we have
\begin{align*}
  \card\;\F_{n-n_0}&= \frac{1}{|R_v^\times|^2}\;
\card\Big(R_v\bs\Big\{(x,y) \in
R_v\times R_v: \begin{array}{c}x R_v+yR_v=R_v\\
  0<|y|_v\leq q_v^{n-n_0}\end{array}\Big\}\Big)
\\&\sim \frac{q^{2\mathfrak g-2}\;q_v^3}{(q-1)^2\,
(q_v^2-1)\,(q_v+1)\;\zeta_K(-1)}
\;q_v^{2n-2n_0}\;.
\end{align*}
For all $n\in\ZZ$ and $[\ga]\in\Ga\!_{\H_\infty}\bs\Ga/
\Ga\!_{\H_\infty}$ outside a finite subset, we have
\[
d(\H_\infty,\ga\H_\infty)>0, \quad m_\ga=1 \quad {\rm and} \quad 
(\varphi^{-1})_*(\Delta_{\Ga\flow{2n}\rho_\ga})=
\Delta_{\Ga r_\ga \Phi^{2n}M}\;.
\]
By Equations \eqref{eq:stabmeasdiscr}, \eqref{eq:expocritnonarch} and
\eqref{eq:psipmstarnonarch}, with $dm$ the counting measure on $\ZZ$,
for every $n_0\in \ZZ$, for $y\in\Ga_H\bs H$ and $m\geq n_0$, we have
\begin{align*}
  d\big((\varphi^{-1})_*\big(\mu^{0+}_{\D^+,2n_0}
  \;_{\mid \;\Ga\bs\gengeod_{\rm even}\XX_v}\big)\big)
(\Theta(y\,\Phi^{-2m})) &=\|\sigma^-_{\D^+}\|\;
d\mu_{\Ga_H\bs H}(y)\,e^{-(\log q_v)\,2m}\,dm\\ &
=\frac{q^{\mathfrak g-1}}{q-1}\;d\mu_{\Ga_H\bs H}(y)\,q_v^{-2m}\,dm\;.
\end{align*}
The end of the proof of Corollary \ref{coro:nonarchFarey} proceeds
now as the one of Corollary \ref{coro:marklof}, replacing Theorem
\ref{theo:quotient} (1) by Theorem \ref{theo:quotient} (3).
\cqfd

{\small \bibliography{../biblio} }

\bigskip
{\small
\noindent \begin{tabular}{l} 
Department of Mathematics and Statistics, P.O. Box 35\\ 
40014 University of Jyv\"askyl\"a, FINLAND.\\
{\it e-mail: jouni.t.parkkonen@jyu.fi}
\end{tabular}
\medskip

\noindent \begin{tabular}{l}
Laboratoire de mathématique d'Orsay, UMR 8628 CNRS,\\
Universit\'e Paris-Saclay,\\
91405 ORSAY Cedex, FRANCE\\
{\it e-mail: frederic.paulin@universite-paris-saclay.fr}
\end{tabular}
}

\end{document}